\documentclass[runningheads]{llncs}
\usepackage[T2A,OT2,OT4,T1]{fontenc}
\usepackage[utf8]{inputenc}
\usepackage[polish,greek,english]{babel}
\RequirePackage{amsfonts}
\usepackage{comment}
\usepackage{multirow}
\usepackage{mathtools}

\DeclarePairedDelimiter\floor{\lfloor}{\rfloor}
\usepackage{appendix}
\usepackage{listings}
\usepackage{xcolor}
\definecolor{mygreen}{RGB}{28,172,0} 
\definecolor{mylilas}{RGB}{170,55,241}
\usepackage{filecontents}
\usepackage[numbers,sort&compress]{natbib}
\usepackage{amsmath}

\usepackage{bbm}
\usepackage{dsfont}
\usepackage{mathtools}

\usepackage{graphicx}
\usepackage{subcaption}

\graphicspath{{./Figures/},{../Figures/},{./Pictures/}}
\usepackage{pgf}
\usepackage{tikz}\usetikzlibrary{arrows,automata}
\usepackage{float}
\usepackage{enumerate}
\usepackage{hyperref}
\hypersetup{
     pdftitle={A measure of~the~importance of~roads based on~topography and traffic intensity}, 
     pdfauthor=Krzysztof Szajowski and Kinga Włodarczyk, 
     colorlinks,
     urlcolor=blue,
     filecolor=magenta,
     citecolor=green, 
     linkbordercolor={1 1 1}, 
     citebordercolor={1 1 1},  
     urlbordercolor={ 1 1 1}  
} 
\RequirePackage[hyperpageref]{backref} 
    \renewcommand*{\backref}[1]{}  
    \renewcommand*{\backrefalt}[4]{
       \ifcase #1 
          No cited.
       \or
          Cited on p. #2.
       \else
          Cited on pp. #2.
       \fi}

\newcommand{\authorcontributions}[1]{%
\vspace{6pt}\noindent{\fontsize{9}{9}\selectfont\textbf{Author Contributions:} {#1}\par}}

\newcommand{\funding}[1]{
\vspace{6pt}\noindent{\fontsize{9}{9}\selectfont\textbf{Funding:} {#1}\par}}

\newcommand{\conflictsofinterest}[1]{%
\vspace{6pt}\noindent{\fontsize{9}{9}\selectfont\textbf{Conflicts of Interest:} {#1}\par}}

\newcommand{\abbreviations}[1]{%
\vspace{12pt}\noindent{\selectfont\textbf{Abbreviations}\par\vspace{6pt}\noindent {\fontsize{9}{9}\selectfont #1}\par}}

\usepackage{glossaries}

\makeglossaries

\newglossaryentry{$Q(t)$}{
name=Survival Function, 
description={It is the reliability function (alternative names: the survivor function): $Q(t)=\bar{F}(t)=1-F(t)$}
}
\newglossaryentry{$phi$}{
name=Structure Function, 
description={It is the structure function: $phi$}
}
\def\subclass#1{\par\addvspace\medskipamount{\rightskip=0pt plus1cm
\def\and{\ifhmode\unskip\nobreak\fi\ $\cdot$
}\noindent\subclassname\ignorespaces #1 \par}}
\DeclareMathSymbol{\preccurlyeq}  {\mathrel}{AMSa}{"34}

\newcommand{\orcid}[1]{\href{https://orcid.org/#1}{\protect\includegraphics[scale=.05]{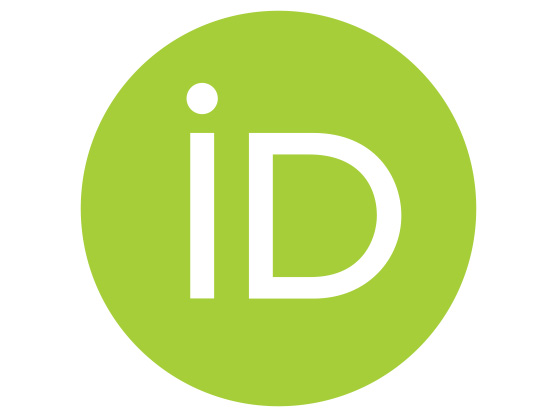}}}
\title{A measure of~the~importance of~roads based on~topography and traffic intensity\thanks{Presented at \textit{Bernoulli-IMS One World Symposium 2020, August 24-28, 2020}. Recorded presentation is \href{https://www.youtube.com/watch?v=joQvFS-r5Kg}{here}.}
}
\author{Krzysztof J. Szajowski\inst{1}\orcid{ 0000-0001-9834-9929} \and
Kinga Włodarczyk\inst{2}\orcid{1111-2222-3333-4444} 
}
\titlerunning{A measure of~the~importance of~roads}
\authorrunning{K.~Szajowski and K. Włodarczyk}
\institute{Wrocław University of Science and Technology, Faculty of Pure and Applied Mathematics, Wybrzeże Wyspiańskiego 27, 50-370 Wrocław, Poland
\email{Krzysztof.Szajowski@pwr.edu.pl}\\
\url{http://szajowski.wordpress.com/} \and
Wrocław University of Science and Technology\\
\email{kingawlodarczyk96@gmail.com}
}
\begin{document}

\maketitle              
\vspace{-2ex}
\begin{abstract}
Mathematical models of~street traffic allowing assessment of~the~importance of~their individual segments for the~functionality of the~stre\-et system is~considering. Based on~methods of~cooperative games and the~reliability theory the~suitable measure is~constructed. The~main goal is~to~analyze methods for assessing the~importance (rank) of~road fragments, including their functions. A~relevance of~these elements for effective accessibility for the~entire system will be~considered.


\subclass{ MSC  68Q80 \and (90B20; 90D80) }

\end{abstract}
\section{\label{intro}Introduction.}
\subsection{\label{IntroHistory}Historical remarks and motivations.} The function of a road network is to facilitate movement from one area to another. As such, it has an important role to play in the urban environment to facilitate mobility. It furthermore determines the accessibility of an (urban) area (together with public transport options). In many studies on the design and maintenance of roads, the authors raise the problem of alternative connections needed to ensure efficient transport between strategic places (cf. \citeauthor{Lin2010:Path}~(\citeyear{Lin2010:Path}), \citeauthor{TacMerMan2012:Hazards}~(\citeyear{TacMerMan2012:Hazards})). It is known that individual segments of the road structure are exposed to various types of threats, resulting in temporary disconnection of such couplings. As a result, the road network determines the quality of life in the analyzed area. Therefore, it is worth trying to define measurable parameters, the quality of road connections, the road system constituting the infrastructure used for transport. Further considerations focus on road systems for road transport. However, the proposed approach can be successfully applied to other similar structures.

When designing, it is worth conducting an analysis of the effects of excluding individual segments and determining the measures that allow for the identification of critical ones. However, the difficulty with this kind of economic appraisal is first of all that it is not easy to measure the valuation of travel time. Different people and organizations value travel time in different ways, depending on many factors such as income, goal of the trip, social background, etc (cf. \citeauthor{Che1981:TTS}~(\citeyear{Che1981:TTS})). It is relatively easier to measure the value travel time than the highway security measure (v. \citeauthor{Sha2012:Secirity}~(\citeyear{Sha2012:Secirity})). The purpose of~the~work is~to~determine the~importance of~roads segment in~road traffic. Consideration will be~commonly known measures of~significance used to~evaluate components of~binary systems. Road topography is~a~long-term process that cannot be~changed in~a~short time. Therefore, it~is~important to~ensure safe road traffic when planning communication infrastructure. To~this end, it~is~important to~introduce objective methods for assessing weak links in~the~road system. Using methods of~stochastic processes and game theory, a~quantitative approach to~the~importance of~various elements of~infrastructure will be~proposed. The introduced connection assessment proposals will be illustrated using information about the actual local road network in the selected city (see Example~\ref{IntroExample1}). 

\subsection{\label{IntroExample1}A motivating example.}	
\begin{figure}[!h!]
{\noindent
\small
\begin{subfigure}[b]{0.5\textwidth}
\begin{center}
\includegraphics[width=0.95\textwidth]{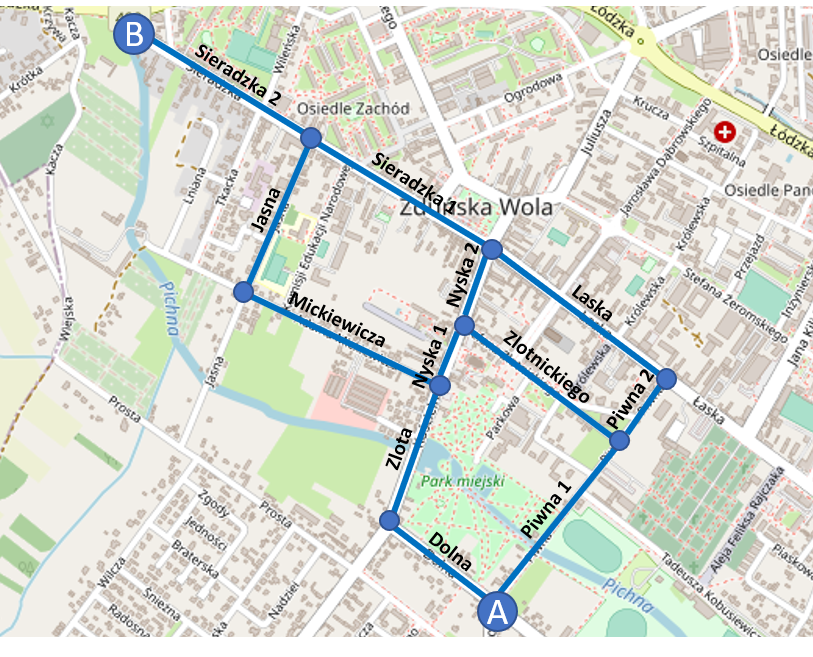}
\caption{\label{fig:main_structureA}Road segments from A to B.}
\end{center}
\end{subfigure}
\hfill\begin{subfigure}[b]{0.5\textwidth}
\centering {\includegraphics[width=0.95\textwidth]{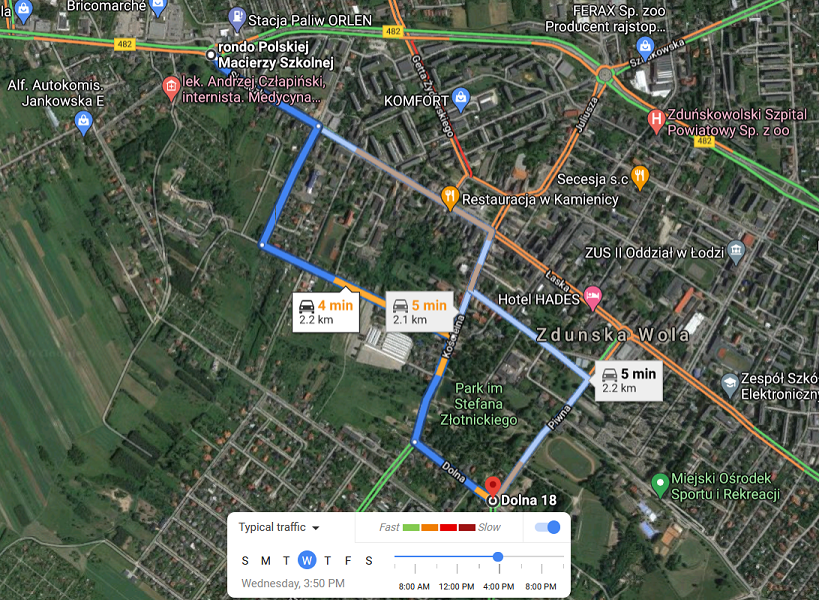}}
\caption{\label{fig:main_structureB} Google Maps presentation of the typical traffic load on the streets.\footnotemark}
\end{subfigure}
}\\[-2ex]

\centering {\caption{\label{fig:main_structure} Analysed traffic network.} }
\end{figure}

In~the~presented work, the~network of~streets ensuring access from point $A$ to~point $B$ in Zduńska Wola will be~treated as~a~system. The~diagram of~the~streets analyzed can be~seen in~Figure \ref{fig:main_structureA}. Let us emphasize that the purpose of modeling is not to reflect the current traffic on the network, as shown in Figure \ref{fig:main_structureB}, but to establish the importance of network elements due to their objective importance for the functioning of the road system.
\footnotetext{Source: Google Maps.}


	In research, there are many measures that allow assessing the~importance of~individual components, based on~the~system structures, lifetimes and reliability of~individual components or~methods of~estimating significance based on~the~methods of~turning on~and off. The~most classic methods based on~reliability theory will be~used in~the~paper. To~this end, the~street network will be~presented in~the~form of~a~system, where each road is~presented as~a~separate component. Then, based on~the~construction of~the~system, the~structure function will be~determined, thanks to~which it~will be~possible to~calculate the~meaning of~individual components and the~corresponding streets. The~next stage will be~defining the~theory of~traffic in~the~context of~significance measures. This area will be~examined in~relation to~the~satisfaction and comfort of~drivers. Drivers satisfaction means that the~system works properly, if~not the~system is~failed. So~as reliability of~particular road we consider probability of~driver's satisfaction from the~journey. A~road connection system in~a~given area should allow transport in~a~predictable time between different points. Extending this time has a~negative effect on~drivers. As a~result, their right ride quality is~compromised and they are more likely to~fail to~comply with the~rules. Therefore, providing drivers with driving comfort and satisfaction is~also important for general road safety. The~use of~this approach can, therefore, be~a~guide for both drivers and road builders planning road infrastructure.

\subsection{\label{IntroPGuide}The paper organization.} The purpose of the research presented here is to implement of various importance measures introduced in the reliability theory to analyze the impact of elements of road networks. The theory related to~significance measures and their use in~traffic theory is~described in~the section~\ref{chap:Importance measures}. There is~a~close relationship between road delays and the~construction and function of~both the~road and the~intersection that forms part of~it. For this purpose, simulations of~vehicle traffic on~the~analyzed roads were performed. The~theory related to~the~method of~modeling vehicle traffic and their behavior at~intersections is~described in~section~\ref{chap:trafic_modelling}. Section~\ref{chap:simulations} describes the~real traffic network, its transfer to~the~simulation model, and the~results obtained in~this way. Then, on~this basis, the~importance of~individual fragments was calculated depending on~the~intensity of~traffic on~these roads. 

Considering this work is~a~look at~the~impact on~the~comfort of~communication of~the~road structure in~connection with traffic without directly referring to~the~behavior of~drivers, which was devoted 
the paper \citeauthor{SzaWlo2020:Divers}~(\citeyear{SzaWlo2020:Divers}). In~these previous works, significant dependence on~traffic quality on~drivers' compliance with applicable rules was shown. Here, a~similar approach was applied to~the~condition of~changing behavior to~incorrect, which may further result in~a~deterioration in~traffic quality. Therefore, the~results obtained show important elements of~the~road network that have an impact on~road safety and properly functioning.  
	\section{ \label{chap:Importance measures}Importance measure. }
	
	The operation of~most systems depends on~the~functioning of~its individual components. It is~important to~ensure the~proper running of~the~entire system. To~this end, it~is~important to~assess the~contribution of~individual components. In road networks, network curves model road segments, intersections, and special places on the road that have a significant impact on the flow of traffic, such as railway crossings, tunnels, bridges, viaducts or road narrowing. In~order to~estimate the~importance of~particular elements, the~concept of~\textit{importance measures} was introduced (for detailed description of the concept and its extension to multilevel elements and systems see review paper by \citeauthor{Overview_importance}~(\citeyear{Overview_importance})). Since 1969~researchers offer various numerical representations to~determine which components are the~most significant for system reliability. It is~obvious that the~greater are these values, the~more this element have on~the~functioning of~the~entire system. The~significance of~individual elements depends on~the~system structure as~well as~the~specificity and failure rate of~individual elements. There are three basic classes of~measures of~importance (v. \citeauthor{Overview_importance}~(\citeyear{Overview_importance}), \citeauthor{Birnbaum1968}~(\citeyear{Birnbaum1968}), \citeauthor{Sre2020:Importance}~(\citeyear{Sre2020:Importance})):
	\begin{enumerate}[i]\itemsep1pt \parskip1pt \parsep1pt
		\item \textbf{Reliability importance measures} subordinate changes in~the~reliability of~the~system depending on~the~change in~the~reliability of~individual elements over a~given period of~time and in~depend on~the~structure of~the~system. 
		\item \textbf{Structural importance measures} are using when just the~structure of~the~system is~known. Depending on~the~position of~the~components in~the~system, their relative importance is~measured. 
		\item\label{TDL} \textbf{Lifetime importance measures} focus on~both, components position in~the~system and lifetime distribution of~each element. According to~\citeauthor{Kuo_Zhu} if~it~is~a~function of~the~time it~can be~classified as~Time-Depend Lifetime(TDL) importance and if~it~is~not a~function of~time we have Time Independent Lifetime(TIL) importance. 
	\end{enumerate}
	
	Moreover, depending on~the~number of~states, systems can be~divided into two types:
	\begin{enumerate}[i]\itemsep1pt \parskip1pt \parsep1pt
		\item \textbf{Binary systems} --- comprised of~$n$ components, where each of~them can have precisely one of~two states. State $\mathbf{0}$ when the~component is~damaged and state $\mathbf{1}$ when is~working. 
		\item \textbf{Multistate systems (MSS)} --- comprised of~$n$ components, which can undergo a~partial failure, but they do not cease to~perform their functions and do not cause damage to~the~entire system. 
	\end{enumerate}
	
	\subsection{\label{sec:importance_in_games}Concepts of~importance measures.} Establishing the hierarchy of components of a complex system has been reduced to measuring the influence of the element state on the status of the entire system. The concept of an element (system) state depends on the context. For the needs of the road network analysis, we assume, similarly to the reliability theory, a binary description of both elements and the system (v. also \citeauthor{Ramamurthy1990}~(\citeyear{Ramamurthy1990})). For this purpose, we will use the known results on significance measures obtained in research on this subject developed in recent years. Importance measures have been developed in~many directions and under many definitions. However, one of~the~most popular areas of~development and application are:
	\begin{itemize}\itemsep1pt \parskip1pt \parsep1pt
		\item theory of~cooperative games (in~simple game)
		\item reliability theory (in~coherent and semi-coherent structure)
	\end{itemize}
	Many methods have been developed to~combine and standardize the~terminology associated with both applications. Therefore, to~begin with, we must briefly mention the~relationship between importance measure theory in~the~context of~both concepts. The~first attempt to~define it~was made by~\citeauthor{Ramamurthy1990}~(\citeyear{Ramamurthy1990}), so~the~following notation was proposed:
	\begin{enumerate}\itemsep1pt \parskip1pt \parsep1pt
		\item[(1)] $\emptyset\in P$, where $P$ is~set of~subset of~$N$;
		\item[(2)] $N\in P$, where $N$ is~finite, nonempty set;
		\item[(3)] $S\subseteq T\subseteq N$ and $S\in P$ imply $T\in P$. 
	\end{enumerate}
	The concepts related to~game theory and reliability theory were compared with each other and on~this basis, it~was possible to~define the~relationship between these concepts. To~begin with, it~is~easy to~see the~relationship between players and components. According to~game theory, we have a~set of~players $N = \{1,2,3,\ldots,n\}$ and a~family of~coalitions $2^N$. In~the~theory of~reliability, we have a~set of~components $N=\{1,2,3,\ldots,n\}$, where the~components and the~entire system can be~in~two states, state 1 for functioning and state 0 for failed. Similarly is~in~game theory, where $\lambda: 2^N \rightarrow \{0,1\}$, which is~applied in~simple game if~on~set $N$ form of~characteristic function fulfils
	\begin{enumerate}\itemsep1pt \parskip1pt \parsep1pt
		\item[(1)] $\lambda(\emptyset)=0$;
		\item[(2)] $\lambda(N)=1$;
		\item[(3)] $S\subseteq T\subseteq N$ implies $\lambda(S)\leq \lambda(T)$.
	\end{enumerate} 
	Here this characteristic function has its counterpart as~a~structure function, and simple games as~semi-coherent structures. In~addition, also winning and blocking coalitions are comparable to~path and cut sets.
	
	In this paper, the~reliability concept of~application of~Importance measures will be~considered, and the~traffic network will be~shown as~a~system. That way, in~the~rest of~this paper we will use reliability terminology. 	
\subsection{\label{sec:ImpMeasRS}Important measures on binary reliability systems.} In the classical approach the system and their elements are binary (v. \citeauthor{Birnbaum1968}(\citeyear{Birnbaum1968}), \citeauthor{Birnbaum1961}(\citeyear{Birnbaum1961})). 
Let the~system comprised of~$n$ components can be~denoted by $	c = (c_1, c_2, ..., c_n)$. The description of the \textit{vector of~component states} (in~the~short \textit{state vector}) $	x = (x_1, x_2, ..., x_n)$, where each $x_i=\chi_W(c_i)$, $c_i\in\{W,F\}$ ($W$ - means the element is functioning; $F$ - means the element is~failed). For state vector, we can use below notations \cite{Birnbaum1968}
	\begin{align*}
	& \vec{x} \leq \vec{y} \mbox{\quad if~\quad} x_i \leq y_i \mbox{\quad for \quad $i \in\{ 1,\ldots,n\}$} \\
	& \vec{x} = \vec{y} \mbox{\quad if~\quad} x_i = y_i \mbox{\quad for \quad $\forall_{i \in \{1,\ldots,n\}}$} \\
	& \vec{x} < \vec{y} \mbox{\quad if~\quad} \vec{x} \leq \vec{y} \mbox{\quad \, and \quad $x \neq y$} \\
	& (1_i,x) = (x_1,x_2,x_3,\ldots,x_{i-1},1,x_{i+1},...,x_n)=(1,x_{-i}) \\
	& (0_i,x) = (x_1,x_2,x_3,\ldots,x_{i-1},0,x_{i+1},...,x_n)=(0,x_{-i}) \\
	& \vec{0} = (0,0,\ldots,0)\mbox{\quad $\vec{1} = (1,1,\ldots,1)$.}
	\end{align*}
	
If~the~structure of~the~system is~known, we can define the~state of~the~system $\phi(\vec{x})$  as~Boolean~function (\textit{structure function} ) of~the~state vector. 
	
	If from $x_i \leq y_i $ for $i \in\{ 1,\dots, n\}$ results $\phi(\vec{x}) \leq \phi(\vec{y}) $, and $\phi(\vec{1}) = 1$,  $\phi(\vec{0}) = 0$,  then we call the~system \textit{coherent}. It is known (v. \citeauthor{Birnbaum1968}~(\citeyear{Birnbaum1968})) that for every $ i = 1,2,\ldots,n $ structure function can be~decomposed as follows:
\begin{equation}\label{eq:structfun2}
	 \phi (\vec{x}) = x_i \cdot \delta_i(\vec{x}) + \mu_i(\vec{x}),    
\end{equation}
where $\delta_i(\vec{x}) = \phi(1_i,\vec{x}) - \phi(0_i,\vec{x})$,  $\mu_i(\vec{x}) = \phi(0_i,\vec{x})$ are independent of~the~state $x_i$ of~the~component $c_i$.	
	
In addition, we can observe situations where the~system can be~functioning even if~some components are failed. The~smallest set of~functioning elements that ensures the~operation of~the~entire system is~called \textit{minimal path}. The~opposite situation is~observed in~the~case of~\textit{minimal cut set}, which is~the~minimum set of~components whose failure cause the~whole system to~fail. We can define the~structure function as~a~parallel structure of~minimal paths. According to~the~definition, this structure is~damaged if, and only if~all of~the~components are failed. So~the~system consists of~$n$ minimal paths series, denoted by~$\rho_i(\cdot)$, for $i = 1,2,\ldots, n$, can be~presented as: 
	\begin{equation}\label{eq:min_path}
	\phi(\vec{x}) = \coprod_{i=1}^{n} \rho_i(\vec{x}) =  1 - \prod_{i=1}^{n} \big[ 1-\rho_i(\vec{x}) \big]. 
	\end{equation}
Similarly, the~structure function can be~presented as~series of~minimal cut sets. So~for $n$ minimal cut parallel structures, marked by~$\kappa_i(\cdot)$, for $i = 1,2,\ldots,n$, structure function looks as~follows:
	\begin{equation}\label{eq:min_cut}
	\phi(\vec{x}) =\prod_{i=1}^{n} \kappa_i(\vec{x}). 
	\end{equation}
	If we simply replace the~minimum paths and minimum cut sets with components, the~formulas \eqref{eq:min_path} and \eqref{eq:min_cut} apply for serial and parallel components.
	
	In most of~the~considerations about the~functioning of~systems, it~is~assumed that elements work independently. Then the~state of~$i$-th element is~a~binary random variable $X_i$ and the~\textit{reliability} that the~element $i$ is~unimpaired will be~denoted by~$p_i$, where
	\begin{equation} \label{eq:probabilities}
	p_i = P(X_i = 1) = 1 - P(X_i = 0).
	\end{equation} 
	We also define the~vector of~reliabilities for $n$ elements by~
	\begin{equation}
	\vec{p} = (p_1,p_2,\ldots,p_n). \label{eq:prob_vec}
	\end{equation} 
	Based on~reliabilities vector and structure function we can define the~probability of~the~system functioning
	\begin{equation}
	P(\phi(x) = 1|\vec{p}) = E[\phi(x)|\vec{p}] = h_{\phi}(\vec{p}). \label{eq:reliability_fun}
	\end{equation}
	For the~structure $\phi(\vec{x})$ function $h_{\phi}(\vec{x})$ is~called \textit{reliability function}.
		
\subsection{\label{sec:reliability}Reliability importance measure. } 	As was introduce, reliability importance measures are based on~changes in~reliabilities of~components and on~the~system structure. This measure first time was introduced by~\citeauthor{Birnbaum1968}~(\citeyear{Birnbaum1968}). At the~beginning, from formulas \eqref{eq:probabilities}, \eqref{eq:prob_vec} and \eqref{eq:structfun2} he express the~reliability function by
$	h_{\phi}(\vec{p}) = p_i \cdot E[\delta_i(X)] + E[\mu_i(X)]$,	where, for every $i = 1,2,\ldots,n$ and according to~equation \eqref{eq:structfun2}, we have
	\begin{equation*}
	\frac{\partial h_{\phi}(\vec{p})}{\partial p_i} = E[\delta_i(\vec{X})] = E\left[\frac{\partial \phi(\vec{X})}{\partial X_i}\right].
	\end{equation*}
	According to~\citeauthor{Birnbaum1968}~(\citeyear{Birnbaum1968}) 
the \textit{reliability importance} of~the~component $c_i$ for structure $\phi(\cdot)$ is defined as $I_{i}(\phi;p) = I_{i}(\phi,1;p) + I_{i}(\phi,0;p)$, where $	I_{i}(\phi,1;\vec{p}) = P\{\phi(X) = 1 | X_i =1;\vec{p}\} - P\{\phi(\vec{X}) = 1;\vec{p}\}$, and $	I_{i}(\phi,0;p) = P\{\phi(X) = 0 | X_i =0;p\} - P\{\phi(X) = 0;p\}$\footnote{$	I_{i}(\phi,1;\vec{p})$ and $	I_{i}(\phi,0;\vec{p})$ are the \textit{reliability importance} of~the~$c_i$ component for functioning and failure of the structure, respectively.}. We have the~following useful identities
	\begin{align*}
	& I_{i}(\phi;1;\vec{p}) = (1-p_i)\cdot\frac{\partial h(\vec{p})}{\partial p_i} = E[(1-X_i)\delta_i(X)] \\
	& I_{i}(\phi;0;\vec{p}) = p_i\cdot\frac{\partial h(\vec{p})}{\partial p_i} = E[X_i\delta_i(X)]\\
	& I_{i}(\phi;\vec{p}) = \frac{\partial h(\vec{p})}{\partial p_i} = E[\delta_i(X)].
	\end{align*}
	
	The~Birnbaum importance measures for $i=1,2,\ldots,n$ have forms (the symbol $\phi$ is droped for short)
	\begin{equation}\label{eq:Birnbaum_Imp1}
	B(i|\vec{p}) = \frac{\partial h(\vec{p})}{\partial p_i} = \frac{\partial [1-h(\vec{p})]}{\partial [1-p_i]}, 
	\end{equation} 
	here $B(i|p)$ is~$p$ dependent. In~case when reliabilities vector $\vec{p}$ is~unknown, we have to~consider structural importance defined for $i=1,2,...,n$ in~the~following way
	\begin{equation}\label{eq:Birnbaum_Imp2}
	B(i) = I_{i}(\phi) = \frac{\partial h(\vec{p})}{\partial p_i}\Bigg\rvert _{p_1=...=p_n=\frac{1}{2}},
	\end{equation}
	this information will be~useful in~the~next section.
		
	\subsection{\label{sec:structural_importance}Structural importance measures.} When we looking for \textit{relevant component} $c_i$ \textit{for the~structure} $\phi(\cdot)$ and the~state vector $\vec{x}$~is~known, we are going as~follow definition $	\delta_i(\vec{x}) = \phi(1_i,\vec{x}) - \phi(0_i,\vec{x}) = 1$.	We can also highlight definitions if~the~component $c_i$ is~\textit{relevant for the~functioning of~structure} $\phi(\cdot)$ at~the~state vector $\vec{x}$ if $	(1-x_i)\cdot\delta_i(\vec{x}) = 1$,	and, if~the~component $c_i$ is~\textit{relevant for the~failure of~structure} $\phi(\cdot)$ at~the~state vector $\vec{x}$ gives $x_i\cdot\delta_i(x) = 1$.	Distinctly, depends on~if the~coordinate $x_i$ of~the~vertex $\vec{x}$ is~equal to~$0$ or~$1$, then $c_i$~is~relevant for functioning or~failure of~the~system.
	
\citeauthor{Birnbaum1968}~(\citeyear{Birnbaum1968}) defined \textit{structural importance measure} of~component $c_i$ \textit{for the~functioning of~the~structure} $\phi(\cdot)$ as $I_{i}(\phi,1) = 2^{-n}\sum_{(x)} (1-x_i) \cdot \delta_i(x)$,	where sum extends on~all combinations $2 ^ n$ of~vertices of~the~state vectors. In~the~similar way is~defined \textit{structural importance measure} of~the~component $c_i$ \textit{for the~failure of~structure} $\phi(\cdot)$ by~$I_{i}(\phi,0) = 2^{-n}\sum_{(x)} x_i \cdot \delta_i(x)$. Finally, by~summarizing, the~\textit{structural importance measure} of~the~component $c_i$~for the~structure $\phi(\cdot)$ is~$I_{i}(\phi) = I_{i}(\phi,1) + I_{i}(\phi,0) = 2^{-n}\sum_{(x)} \delta_i(x)$.
	
	
	\citeauthor{barlow}~(\citeyear{barlow}) used a~more extended approach to~structural measures. Their point of~view assumes that all components have a~continuous lifetime distribution, denoted by~$ F_i $, for $ i = 1,2, \ldots, n $. It is~possible to~calculate the~probability of~a~system failure caused by~the~$ c_i $ component. For the~$ c_i $ component, which is~described by~the~distribution $ F_i $ and the~density function $ f_i $, the~probability that a~system failure at~time $t$~was caused by~the~$ c_i $ component can be~described as~follows
	\begin{equation}\label{eq:lifetime}
	\frac{[h(1_i,\bar{F}(t)) - h(0_i,\bar{F}(t))] f_i(t)}{\sum_{k=1}^{n}[h(1_k,\bar{F}(t)) - h(0_k,\bar{F}(t))]f_k(t)}.
	\end{equation}
	In the~consequence of~\eqref{eq:lifetime}, it~obvious to~define the~probability that failure of~the~system in~$[0,t]$ was caused by~the~$c_i$ component is
	\begin{equation*}
	\frac{\int_{0}^{t}[h(1_i,\bar{F}(u)) - h(0_i,\bar{F}(u))] d F_i(u)}{\int_{0}^{t} \sum_{k=1}^{n}[h(1_k,\bar{F}(u)) - h(0_k,\bar{F}(u))]d F_k(u)}.
	\end{equation*}
	Here, if~$t \to \infty$, then we obtain that the~system finally failed it~was caused by~the~component~$c_i$. In~this case, we have to~note that the~denominator is~equal to~1. This limit is~taken as~a~definition of~\emph{component importance}.
	
Importance measures according Barlow and Proschan definition we will denoted by~$I_{i}^{BS}(\phi)$. We have
	\begin{equation}\label{eq:BP_importance}
	I_{i}^{BP}(\phi) =  \int\limits_{0}^{1} [h(1_i,p) - h(0_i,p)] d p,
	\end{equation}
	where $(1_i,p)$ and $(0_i,p)$ is~a~probability vector where $i$-th component has probability equal 1 or~0, relatively.
	
	For further calculations, let us remind quick note from Section \ref{sec:ImpMeasRS}, that \textit{minimal path} is~the~minimal set of~elements, which ensures the~proper functioning of~the~system. 
	Based on~this we can define  \textit{critical path set} for component $c_i$ as $\{i\} \cup \{j | x_j = 1, i \neq j\}$. In this way, information about the~system is~functioning or~failed is~determined by~the~$c_i$~component functions or~fails. A~\textit{critical path vector} (or \textit{set}) for the~component~$c_i$, and its size $r$, we have $1 + \sum_{i \neq j} x_j= r$, for $r = 1,2,\ldots,n$. The~formula for counting the~number of~vectors of~critical paths for the~component $c_i$ with size $r$ is~the~following 
	\begin{equation*}
	n_r(i) =  \sum_{\sum_{i\neq j}x_j=r-1} [ \phi(1_j,x) - \phi(0_j,x) ].
	\end{equation*}
	Finally, we can define the~structural importance of~the~component $c_i$ using the~number of~vectors of~critical paths $n_r(i)$ as~follows 
	\begin{equation}\label{eq:BP_importance2}
	I_{i}^{BP}(\phi) = \sum_{r=1}^{n} n_r(i) \cdot \frac{(r-1)! (n-r)!}{n!}.
	\end{equation}
	The equation \eqref{eq:BP_importance2} can be~also presented in~two more interesting expressions. The~first expression is~the~following 
	\begin{equation*}
	I_{i}^{BP}(\phi) = \frac{1}{n}  \sum_{r=1}^{n} n_r(i) \tbinom{n-1}{r-1}^{-1},
	\end{equation*}
	where $n_r(i)$ describes the~number of~vectors of~critical paths with size $r$. The~denominator in~the~above equation represents the~amount of~results in~which precisely $r-1$ components are in~operation among the~$n-1$ components without $c_i$ component. 
	Second additional representation of~equation \eqref{eq:BP_importance2} can be~written as~follows
	\begin{equation*}
	I_{i}^{BP}(\phi) = \int\limits_{0}^{1} \Big[ \sum_{r=1}^{n} n_r(i) \cdot \tbinom{n-1}{r-1}^{-1} \tbinom{n-1}{r-1} \cdot (1-p)^{n-r} \cdot p^{r-1}\Big]dp,
	\end{equation*}
	here $\tbinom{n-1}{r-1} (1-p)^{n-r} p^{r-1}$ means the~probability that from the~$n-1$ components without $c_i$~component, $r-1$ elements are functioning. What's more, $n_r(i) \tbinom{n-1}{r-1}^{-1}$ means the~probability that $r-1$ functioning elements including $c_i$ component determine the~critical path set to~the~$c_i$ component. So~multiplication of~them means the~probability that components $c_i$~is~responsible for system failure and integral of~it over $p$ is~that reliability for the~component $c_i$ is~a~uniform distribution  $p \sim \mathcal{U}(0,1)$.
	
	As it~was written at~the~beginning in~Section \ref{sec:importance_in_games} there is~a~big connection between the~concepts related to~game theory and the~theory of~reliability. The~measure introduced by~Barlow and Proschan is~an example of~this. This definition is~reflected in~cooperative games as~Shapley's value, which informs what profit a~given coalition player can expect, taking into account his contribution to~any coalition. 
	
\subsection{\label{sec:traffic_importance}Importance measures of road segments based on~traffic flow in example~\ref{IntroExample1}.} 
\begin{figure}[tbh!]
		\centering
		\includegraphics[width=12cm]{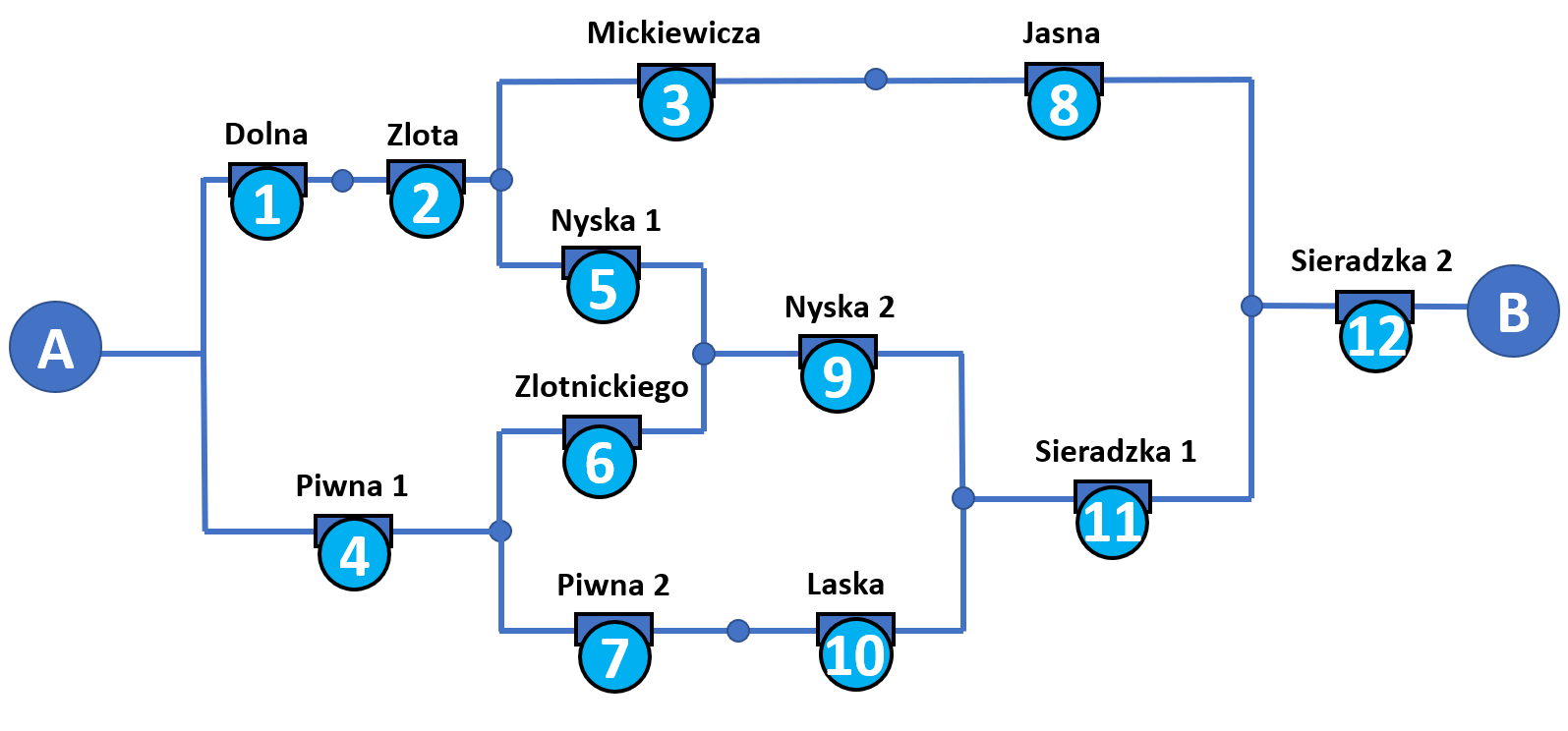} 
		\caption{\label{fig:system}Analysed traffic network presented in~system form.} 
	\end{figure}

As was said in~section~\ref{intro}, binary systems are considered. The~analyzed system is~a~street network allowing access from $A$ to~$B$, it~is~possible in~several ways. We assume that drivers drive only from $A$ to~$B$, straight, without unnecessary U-turns on~the~route. Streets were presented at~the~beginning in~Fig.~\ref{fig:main_structureA} and can be~transform to~the~form of~the~system (a scheme) as on~Fig.~\ref{fig:system}.
	
Based on~the~system representation of~the~streets network we can determine the~structure function. As we know, the~structure function can be~defined using either \textit{minimal path set} or~\textit{minimal cut set}. So~for the~given structure both sets are presented in~the tables \ref{tab:minpathset} and \ref{tab:mincutset}.	
\hspace{-2ex}
	\begin{table}[tbh!]
		\begin{minipage}{.45\linewidth}
			\centering{\small
			\caption{\label{tab:minpathset}Minimal path set.}
			\begin{tabular}{|c|c|} 
				\hline
				Path & Elements    \\ [0.5ex] 
				\hline
				1 & 1 2 3 8 12         \\ 
				2 & 1 2 5 9 11 12        \\
				3 & 4 6 9 11 12        \\
				4 & 4 7 10 11 12          \\
				\hline
			\end{tabular} }
		\end{minipage}
		\begin{minipage}{.45\linewidth}
			\centering{\small
			\caption{\label{tab:mincutset}Minimal cut set.} 			
			\begin{tabular}{|c|c||c|c|} 
				\hline 
				Cut & Elements & Cut & Elements \\ [0.5ex] 
				\hline
				1 & 1 4  &   11 & 3 9 10     \\ 
				2 & 2 4    &  12 & 8 9 7     \\ 
				3 & 1 6 7    &  13 & 8 9 10    \\ 
				4 & 2 6 7    &  14 & 3 4 9    \\ 
				5 & 4 5 3    &  15 & 4 8 9    \\ 
				6 & 4 5 8    &  16 & 3 11   \\
				7 & 1 6 10   &  17 & 8 11   \\
				8 & 2 6 10   &  18 & 1 11  \\
				9 & 3 5 6 7   &  19 & 2 11   \\
				10 & 3 9 7   &   20 & 12   \\
				\hline
				
			\end{tabular}}
		\end{minipage}
	\end{table}
	
Based on~tables \ref{tab:minpathset} and \ref{tab:mincutset}, it~is~possible to~define minimal path series structures represented by~the~following equations
	\begin{align*}
	\rho_1(x) &= \prod_{\{1, 2, 3, 8,12\}}\hspace{-1em}x_i= x_1\cdot x_2\cdot x_3 \cdot x_8 \cdot x_{12}
	\hspace{-.5em} &\rho_2(x) = \prod_{\{1, 2, 5,9,11,12\}}\hspace{-1em}x_i 
	&\quad\\ 
	\rho_3(x) &= \prod_{\{4, 6, 9,11,12\}}\hspace{-1em}x_i 
	\hspace{-.5em} &\rho_4(x) = \prod_{\{4, 7, 10,11,12\}}\hspace{-1em}x_i &\quad\\
	\intertext{and minimal cut parallel structures described as~follows}
\kappa_1(x) &=\coprod_{\{1,4\}}x_i= x_1 \amalg x_4  \hspace{-1.5em} &\kappa_2(x) =\coprod_{\{2,4\}}x_i\quad &\kappa_3(x) =\coprod_{\{1, 6, 7\}}x_i\\
\kappa_4(x) &=\coprod_{\{2, 6, 7\}}\hspace{-.5em}x_i= x_2 \amalg x_6 \amalg x_7 \hspace{-.5em}&\kappa_5(x) =\coprod_{\{4,5,3\}}x_i\ 
& \kappa_6(x) =\coprod_{\{4,5,8\}}\hspace{-.5em}x_i   \\
\kappa_7(x) &=\coprod_{\{1,6,10\}}x_i= x_1 \amalg x_6 \amalg x_{10}\     &\kappa_8(x) =\coprod_{\{2,6,10\}}x_i\ &\kappa_9(x)=\coprod_{\{3,5,6,7\}}x_i\\
\kappa_{10}(x) &=\coprod_{\{3,9,7\}}x_i= x_3 \amalg x_9 \amalg x_7\  &\kappa_{11}(x) =\coprod_{\{3,9,10\}}x_i\ &\kappa_{12}(x) =\coprod_{\{8,9,7\}}x_i \\
	\kappa_{13}(x) &=\coprod_{\{8,9,10\}}x_i= x_8 \amalg x_9 \amalg x_{10}\      &\kappa_{14}(x) =\coprod_{\{3,4,9\}}x_i\ &\kappa_{15}(x) =\coprod_{\{4,8,9\}}x_i\\
	\intertext{}
	\kappa_{16}(x) &=\coprod_{\{3,11\}}x_i= x_3 \amalg x_{11}\ &\kappa_{17}(x) =\coprod_{\{8,11\}}x_i\ &\kappa_{18}(x) =\coprod_{\{1,11\}}x_i \\
	\kappa_{19}(x) &=\coprod_{\{2,11\}}x_i= x_2 \amalg x_{11}\     &\kappa_{20}(x) =\coprod_{\{12\}}x_i= x_{12}\ & 
	\end{align*}
	From the~definition in~equation \eqref{eq:min_path} and based on~the~above equations, we can write the~structure function of~the~presented system as~follows
	\begin{align*}
	\phi(x) &= \rho_1(x) \amalg \rho_2(x) \amalg \rho_3(x) \amalg \rho_4(x) = \\
	& = 1 - (1-\rho_1(x))(1-\rho_2(x))(1-\rho_3(x))(1-\rho_4(x)) 
	\end{align*}
	In addition, our structure function can be~also expressed by~the~series of~minimal cut structures 
	\begin{equation*}
	\phi(x) = \prod_{i=1}^{20} \kappa_i(x).  
	\end{equation*}
	And finally, thanks to~equation \eqref{eq:reliability_fun}, we can write the~reliability function of~the~analyzed system
	\begin{align} \label{eq:my_struc_fun}
	h_{\phi}(p)  &= 1 - (1 - \prod_{\{1, 2, 3,8,12\}}p_i )(1 -\prod_{\{1, 2, 5,9,11,12\}}p_i )(1 - \prod_{\{4, 6, 9,11,12\}}p_i  )\\
	\nonumber&\quad \times(1 - \prod_{\{4, 7, 10,11,12\}}p_i ),
	\end{align} 
	where $p_i$, for $i = 1,2,\ldots,12$, are some probabilities, which definition will be~introduce in~next section.
\subsection{\label{sec:reliabilityImp_traffic}Reliability importance applied to road networks.} To consider reliability importance we need to~define what exactly means that system is~functioning or~failed. We assume that the~condition of~the~system's functioning is~the~comfort and satisfaction of~drivers.  The~state \textbf{1} will mean the~driver's satisfaction with a~given road section or~route, the~state \textbf{0} --- dissatisfaction. For drivers, the~measure of~satisfaction is~the~travel time on~a~given section of~the~road, and more precisely, the~realization of~the~road according to~the~planned travel time. Drivers want to~finish the~journey in~the~shortest possible time. The~excess of~this time, i.e. the~delay on~a~given section of~the~road after exceeding a~certain critical level causes dissatisfaction of~drivers with the~journey. This critical level that causes dissatisfaction may be~different for each driver and is~close to~the~lifetime. Weibull distribution is~often used to~represent the~lifetime of~objects. A~similar approach was used in~a~paper published by~\citeauthor{FanJiaTianYun2014:Weibul} in~\citeyear{FanJiaTianYun2014:Weibul}. The~cited article considered a~situation when, while waiting before entering the~intersection, the~waiting time for a~given driver exceeded a~certain critical value, the~driver stopped complying with traffic rules. Like here, this critical value was determined by~Weibull distribution. The~variable from the~Weibull distribution can be~represented as~the~cumulative distribution function given by~the~following formula:
	\begin{equation*} 
	F(t) = \left\{ \begin{array}{ll}
	1-\exp \left\{ -\left( \frac{t}{\lambda} \right)^k \right\}, &\text{for $t>0$,}\\
	0, & \textrm{otherwise.}\\
	\end{array} \right.
	\label{dystrybuanta}
	\end{equation*}
	Based on~the~cumulative distribution function, it~is~possible to~calculate the~reliability function, i.e. the~function that tells the~probability of~correct functioning of~an object. We parametrize the segments by the acceptable delay time $t$ by the driver. The population of the diver is not homogeneous. The acceptable delay is the random variable with some distribution $\Pi$. The delay of travel $\tau$ is a consequence of various factors. Let us assume that its cumulative distribution is $F(t)$. We will say that the segment is reliable or works for given driver with accepted delay $t$ if $\tau(\omega)\geq t$. The probability $Q(t)$ of the event is the subjective driver reliability of the segment. Its expected probability with respect to $\Pi$, $p=\int_0^\infty Q(u)d\Pi(u)$ is mean reliability of the segment.     
For the homogeneous class of drivers the delay time $\xi$ on~given road section is common value for all drivers, so~(mean) reliability will means probability that for assumed delay time driver is~satisfied from the~journey. Therefore, according to~the~theory of~\emph{importance measures}, $p_i$, which indicates the~reliability of~the~segment will be~determined as~probability of~driver's satisfaction and $1-p_i$ will means probability that driver is~dissatisfied of~journey for assumed delay time. It can be~determine as~following formula:
	\begin{equation*}
	p_i = P(X_i=1 |t = \xi) = 1-P(X_i = 0 |t = \xi) = Q(\xi),
	\end{equation*}
	where $ \xi $ is~the~delay time, $ X = 1 $ means the~driver is~satisfied with the~road, $ X = 0 $ -- he is dissatisfied. The~paper adopts the~same Weibull distribution parameters as~in~paper by~\citeauthor{FanJiaTianYun2014:Weibul}~(\citeyear{FanJiaTianYun2014:Weibul}), i.e. $ \lambda $ = 30, $ k = 2.92 $. 
	When we know the~relationship between street reliability from the~time of~delays, we can define the~reliability of~these route fragments for a~given traffic intensity. Different road sections react differently to~increasing traffic intensity, which is~why reliabilities will be~different. Using simulations we determine the~dependence of~traffic intensity and delay times.
	
	
\subsection{Continuing example of section~\ref{IntroExample1}}	
\begin{figure}[tbh!]
		\centering
		\includegraphics[width=7cm]{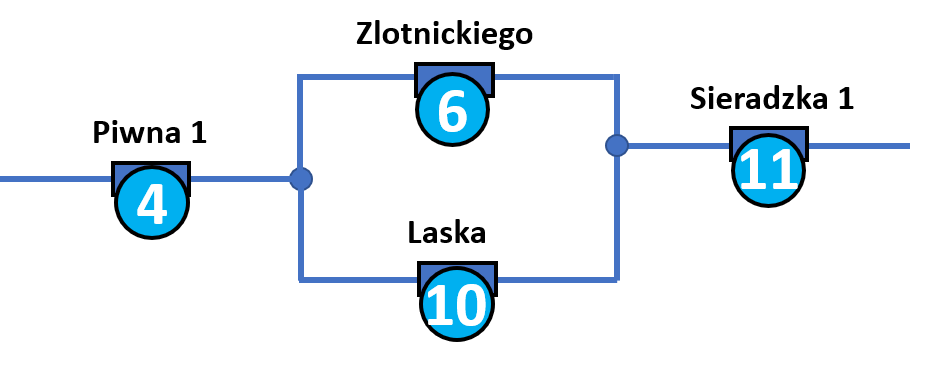} 
		\caption{Short version of~analysed scheme of~traffic network.} \label{fig:short_example}
	\end{figure}

	We will now proceed to~briefly introduce these definitions on~a~simple example. Let us assume that we have a~shortened road network scheme limited to~Piwna~1, Zlotnickiego, Laska, Sieradzka~1 streets. This scheme is~presented in~the~way shown in~Figure~\ref{fig:short_example}. 
	Here we have the~components $c_4$ i $c_{11}$ in~series, and the~components $c_6$ and $c_{10}$ in~parallel. So~we can define this system as~"$k$--out--of--$n$" structure, where $ n $ is~number of~all components, $ k $ means number of~components in~series, and $ n-k $ is~the~number of~components in~parallel. On this basis, we can define the~structure function as $	\phi(\vec{x}) = x_4 \cdot (1-(1-x_6)\cdot (1-x_{10}))\cdot x_{11}$, and the~system reliability function corresponding to~the~above $h(\vec{p}) = p_4 \cdot (1-(1-p_6)\cdot (1-p_{10}))\cdot p_{11}$.
	
	To begin with, we assume that the~reliability of~individual components is~unknown, so~only structural measures of~significance can be~calculated. They will be~calculated based on~the~definitions introduced in~Section \ref{sec:structural_importance}. Two proposals for structural measures have been introduced: first  proposed by~Birnbaum, which assume that each reliability $p_i$~of~components $c_i$, for $i = 1,2,\ldots,n$ are the~same and equal to~$\frac{1}{2}$ and second the~Barlow and Proschan Importance Measures, which is~define for $p \in [0,1]$. So~using this theory and definitions in~equation \eqref{eq:Birnbaum_Imp2} for Birnbaum Importance and in~\eqref{eq:BP_importance} for Barlow and Proschan Importance, we count the~importance of~analyzed components. Obtained results are presented in~Table \ref{tab:importance_struc_theor}.
	\begin{table}[h!]
		\centering{\small
		\caption{Structural importance of~roads in~the~analyzed system.}
		\label{tab:importance_struc_theor}
		\begin{tabular}{|c|c|c|c|}
			\hline 
			\multirow{2}{0.5cm}{\textbf{Id}} & \textbf{Street} & \textbf{Birnbaum} & \textbf{Barlow-Proschan} \\
			& \textbf{name} & \textbf{Importance} $B(i;\phi)$ & \textbf{Importance} $I_i^{BP}(\phi)$ \\ 
			\hline 
			4 & Piwna~1& 0.375 & 0.4167\\ 
			\hline 
			6 & Zlotnickiego & 0.125 & 0.0833 \\ 
			\hline 
			10 & Laska & 0.125 & 0.0833 \\ 
			\hline 
			11 & Sieradzka~1 & 0.375 & 0.4167 \\ 
			\hline 
		\end{tabular} }
	\end{table}
	We see that roads connected in~series are more important than roads connected in~a~parallel way. This is~consistent with the~logic, if~one of~the~parallel roads is~blocked, you can always choose a~different route that will allow you to~reach your destination. For streets in~a~serial connection, this is~not possible. We also note the~differences in~the~values of~the~importance measures calculated using the~Birnbaum and Barlow-Proschan definitions, this is~because the~first measure is~calculated for the~constant reliability of~the~elements equal to~$\frac{1}{2}$, so~it~only examines the~relationship between element positions. The~second measure takes into account, apart from the~structure itself, also the~variability of~reliability of~individual elements.
	
	Now we will examine the~reliability importance measures for the~simplified system shown in~Figure \ref{fig:short_example}. Let us assume that for a~given traffic intensity, we have certain delay times, on~this basis, we will calculate the~probability that drivers are still satisfied with the~travel along the~road, i.e. the~reliability of~the~road. Next, for these values, using the~formula \eqref{eq:Birnbaum_Imp1}, we calculate the~value of~the~measures of~significance defined by~Birnbaum. The~assumed delay times, as~well as~the~corresponding reliability and importance, will be~presented in~Table \ref{tab:relabilityImp:example1}.
	\begin{table}[h!]
		\centering{\small
		\caption{Hypothetical calculations of~importance measures in~the~example system.}
		\label{tab:relabilityImp:example1}
		\begin{tabular}{|c|c|c|c|c|}
			\hline 
			\multirow{2}{1cm}{\textbf{Id}} & \textbf{Street} & \textbf{Delay} & \textbf{Probability of} & \textbf{Importance} \\
			& \textbf{name} & $\xi$ & \textbf{satisfaction} $Q(\xi)$ & $B(i|p)$\\ 
			\hline 
			4 & Piwna~1& 25 s & 0.5559 & 0.8513 \\ 
			\hline  
			6 & Zlotnickiego & 20 s & 0.7363  & 0.0025\\ 
			\hline 
			10 & Laska & 5 s & 0.9947  & 0.1249 \\ 
			\hline 
			11 & Sieradzka~1 & 16 s & 0.8526  & 0.5551 \\ 
			\hline  
		\end{tabular} }
	\end{table}
	As we can see, with a~road delay of~about 25 seconds, the~likelihood of~driver satisfaction is~close to~$ \frac{1}{2} $, and in~the~case of~delays of~about 5 seconds, drivers do not experience almost the~negative effects of~a~slowdown in~traffic. With such reliability of~streets and with such a~scheme, it~is~easy to~notice some issues: streets in~a~parallel position have less contribution to~potential nervousness or~driver satisfaction than in~a~serial connection, in~addition, in~the~case of~streets in~a~series, those with less reliability are more important, so~these should be~paid greater attention to~maintain proper traffic quality. In~the~case of~streets in~parallel connection, streets with greater reliability are more important. It~is~logical that drivers knowing that the~road is~a~better way will choose it, so~it~is~important to~constantly maintain it~in~good condition because when it~fails the~whole connection will lose much reliability.

\subsection{\label{struc_traff}Structural importance of~real traffic network.} 	As was presented in~the~previous sections, if~the~reliability of~individual components is~unknown, it~is~possible to~use structural measures of~significance. Therefore, we~will begin our considerations about the~analyzed system by~calculating the~structural significance of~individual roads. In~the~same way as~in~the~previous section, the~definition of~significance measures introduced by~Birnbaum, and Barlow and Proschan presented in section~\ref{sec:reliability} by~the~equations \eqref{eq:Birnbaum_Imp2} and \eqref{eq:BP_importance}, respectively, were used. The~results obtained are presented in~Table \ref{table:importanceresults}.
	\begin{table}[h!]
		\centering{\small
		\caption{Structural importance of~roads in~the~analyzed system.}
		\label{table:importanceresults}
		\begin{tabular}{|c|c|c|c|}
			\hline 
			\multirow{2}{0.5cm}{\textbf{Id}} & \textbf{Street} & \textbf{Birnbaum} & \textbf{Barlow-Proschan} \\
			& \textbf{name} & \textbf{Importance} $B(i;\phi)$ & \textbf{Importance} $I_i^{BP}(\phi)$ \\ 
			\hline 
			1 & Dolna & 0.0861 & 0.0973\\ 
			\hline 
			2 & Zlota & 0.0861 & 0.0973\\ 
			\hline 
			3 & Mickiewicza &  0.0577& 0.0531\\ 
			\hline 
			4 & Piwna~1& 0.1155 & 0.1202\\ 
			\hline 
			5 & Nyska~1& 0.0284 & 0.0338\\ 
			\hline 
			6 & Zlotnickiego & 0.0577 & 0.0531 \\ 
			\hline 
			7 & Piwna~2 & 0.0577 & 0.0531 \\ 
			\hline 
			8 & Jasna & 0.0577 & 0.0531 \\ 
			\hline 
			9 & Nyska~2&  0.0861 & 0.0973 \\ 
			\hline 
			10 & Laska & 0.0577 & 0.0531 \\ 
			\hline 
			11 & Sieradzka~1 & 0.1439 & 0.1882 \\ 
			\hline 
			12 & Sieradzka~2 & 0.2016 & 0.3690\\ 
			\hline 
		\end{tabular}} 
	\end{table}
	\noindent We see that the~results of~both measures are similar. As was expected, the~most important for the~entire route is~street Sieradzka~2, because each route finally leads along this street, for B-P importance for these streets is~bigger than for B-importance. Next, the~most important part of~the~route is~Sieradzka~1, we see that 3 of~4 ways to~obtain point~$B$~are going by~this street. For this street, Birnbaum's value is~smaller than the~Barlow-Proschan's value. The~importance of~Piwna~1 is~the~last value of~importance bigger than $0.1$, anyway similar to~this value are importances of~Dolna, Zlota, and Nyska~2. Surmise that the~significance of~the~Zlota and Dolna will be~close to~the~value calculated for Piwna~1 was not difficult. However, it~is~not so~easy to~guess the~similarity of~the~importance of~Nyska~2 street to~Dolna and Zlota. Mickiewicza, Zlotnickiego, Piwna~2, Jasna, Laska and Nyska~1 streets have the~smallest contribution to~the~proper functioning of~the~entire connections between $A$ and $B$.

	\section{\label{chap:trafic_modelling}Traffic modelling} 	
	\subsection{Review and history.} Traffic modeling is~a~particularly complex issue. There are both modeling of~individual phenomena occurring on~roads and entire road networks. The~first research into vehicle movement and traffic modeling theory began with the~work of~Bruce D. Greenshields(\citeyear{greens}). On the~basis of~photographic measurement methods, he proposed basic and empirical relationships between flow, density, and speed occurring in~vehicle traffic. Next, \citeauthor{LigWhi1955:Kinematic}~(\citeyear{LigWhi1955:Kinematic}) and \citeauthor{Rich1956:Shock}~(\citeyear{Rich1956:Shock}) introduced the~first theory of~movement flow. They presented a~model based on~the~analogy of~vehicles in~traffic and fluid particles. Interest in~this field has increased significantly since the~nineties, mainly due to~the~high development of~road traffic. As a~result, many models were created describing various aspects of~road traffic. As a~result, many models were created describing various aspects of~road traffic and focusing on~different detail models, we can distinguish:
	\begin{itemize}\itemsep1pt \parskip1pt \parsep1pt
		\item microscopic models
		\item mesoscopic models
		\item macroscopic models
	\end{itemize}
	The differences in~the~models are at~the~level of~aggregation of~modeled elements. Mesoscopic models based mainly on~gas kinetic models. Macroscopic models based on~first and second-order differential equations, derived from Lighthill-Whitham-Richards(LWR) theory. Microscopic models focus on~the~simulation of~individual vehicles and their interactions. One of~the~most popular are car-following models and cellular automata models, the~last is~used in~this paper. The~most popular cellular automata traffic model is~the~Nagel-Schereckenberg (\citeyear{nasch}) model , but also very interesting model is~LAI~model (cf.~\citeauthor{LarAlvIca2010:Cellular}(\citeyear{LarAlvIca2010:Cellular})), which is~more advance than NaSch model. LAI~model is~used in~this paper, therefore, in~the~next section theory about cellular automata will be~introduced and later will be~a~more detailed description of~LAI~model.

	\subsection{Cellular automaton}\label{automaty} 
	Janos von Neumann, a~Hungarian scientist working at~Princeton, is~the~creator of~cellular automata theory. In~addition, the~development of~this area was significantly influenced by~the~Lviv mathematician Stanislaw Ulam, who is~responsible for discrediting the~time and space of~automats and is~considered to~be the~creator of~the~definition of~cellular automats as~"imaginary physics"\cite{automaty3}. According to~a~book written by~\citeauthor{ksiazkaautomaty}, cellular automata can reliably reflect many complex phenomena with simple rules and local interactions. Cellular automata are a~network of~identical cells, each of~which can take one specific state, with the~number of~states being arbitrarily large and finite. The~processes of~changing the~state of~the~cells run parallel and according to~the~rules. These rules usually depend on~the~current state of~the~cell or~the~state of~neighboring cells. From the~mathematical point of~view, cellular automatas are defined by~the~following parameters \cite{automaty2} \cite{automaty1}:
	\begin{itemize}\itemsep1pt \parskip1pt \parsep1pt
		\item \textbf{State space} --- a~finite, $k$-element set of~values defined for each individual cell.
		\item \textbf{Cell grid} --- discrete, $ D $-dimensional space divided into identical cells, each of~which at~a~given time $ t_h $ has one, strictly defined state of~all possible $ k $ states. In~the~case of~the~$ 2D $ network, the~cell status at~$ i $, $ j $ is~indicated by~the~symbol $ \sigma_ {ij} $.
		\item \textbf{Neighborhood} --- parameter determining the~states of~the~nearest neighbors of~a~given cell $ ij $, marked with the~symbol $ N_{ij} $.
		\item \textbf{Transition rules} --- rules determining the~cell state in~a~discrete time $ t_{h + 1} $ depending on~the~current state of~this cell and the~states of~neighboring cells. The~state of~the~cell in~the~next step is~presented in~the~following relationship:
		\begin{equation*} 
		\sigma_{ij}(t_{h+1}) = F\left(\sigma_{ij}(t_h),N_{ij}(t_h) \right),
		\end{equation*}
		where:\\
		\begin{tabular}{p{3.7em} l}
			$\sigma_{ij}(t_{h+1})$ & --- cell state in~position $i$,$j$ in~step $t_{h+1}$, \\ 
			$\sigma_{ij}(t_{h})$ & ---  cell state in~position $i$,$j$ in~step $t_{h}$, \\  
			$N_{ij}(t_h)$ & --- cells in~the~neighborhood of~a~cell in~position $i$,$j$ in~step $t_{h}$.  
		\end{tabular}
	\end{itemize}
	The way the~cell neighborhood is~defined has a~significant impact on~the~calculation results. The~most common are two types:
	\begin{itemize}\itemsep1pt \parskip1pt \parsep1pt
		\item \textbf{Von Neumann neighborhood} Each cell is~surrounded by~four neighbors, immediately adjacent to~each side of~the~cell being analyzed. The~neighborhood for $ i, j $ constructed in~this way is~as follows:
		$$ N_{i,j}(t_h) =
		\left( \begin{array}{ccc}
		& \sigma_{i-1,j}(t_h) &  \\[1ex]
		\sigma_{i,j-1}(t_h) & \mathbf{\sigma_{i,j}(t_h)} & \sigma_{i,j+1}(t_h) \\[1ex]
		&  \sigma_{i+1,j}(t_h) & 
		\end{array} \right) $$
		\item \textbf{Moore neighborhood} Each cell is~surrounded by~eight neighbors, four directly adjacent to~the~sides of~the~analyzed cell, and four on~the~corners of~the~analyzed cell. The~neighbor cell matrix for $ i $, $ j $ looks like this:
		$$ N_{i,j}(t_h) =
		\left( \begin{array}{ccc}
		\sigma_{i-1,j-1}(t_h) & \sigma_{i-1,j}(t_h) & \sigma_{i-1,j+1}(t_h) \\[1ex]
		\sigma_{i,j-1}(t_h) & \mathbf{\sigma_{i,j}(t_h)} & \sigma_{i,j+1}(t_h) \\[1ex]
		\sigma_{i+1,j-1}(t_h) &  \sigma_{i+1,j}(t_h) & \sigma_{i+1,j+1}(t_h)
		\end{array} \right) $$
	\end{itemize}
	There are also modifications to~the~above types, such as~the~combined neighborhood of~Moore and von Neumann, as~well as~numerous modifications to~the~Moore neighborhood itself, and a~different way defined by~Margolus to~simulate falling sand.
	
	In addition, boundary conditions are an important aspect of~cellular automata theory. Since it~is~impossible to~produce an infinite cellular automaton, some of~the~simulations would be~impossible because with the~end of~the~automaton's grid the~history of~a~given object or~group of~objects would end. For this purpose, boundary conditions at~the~ends of~the~grid were introduced. There are the~following types of~boundary conditions:
	\begin{itemize}\itemsep1pt \parskip1pt \parsep1pt
		\item periodic boundaries --- cells at~the~edge of~the~grid behind neighbors have cells on~the~opposite side. In~this way, the~continuity of~traffic and ongoing processes is~ensured. 
		\item open boundaries --- elements extending beyond the~boundaries of~the~grid cease to~exist. This is~used when new objects are constantly generated, which prevents too high density of~objects on~the~grid.
		\item reflective boundaries --- on~the~edge of~the~automaton a~border is~created, from which the~simulated objects are reflected, most often it~serves to~imitate the~movement of~particles in~closed rooms.
	\end{itemize}
	In the~next section the~model using cellular automata used in~the~simulation will be~presented. Open boundary conditions are used in~our simulations. After leaving the~street, vehicles disappear. This is~in~line with logic, new vehicles are constantly appearing and disappearing on~the~roads. The~applied neighborhood is~a~modified version of~the~presented neighborhoods, because vehicles as~their neighbors take those vehicles that are nearby, and more specifically the~nearest vehicle on~the~road, even if~it~is~not directly adjacent to~the~analyzed vehicle, and also cars move by~more cell. We can assume that it~is~a~more extended version of~the~Von Neumann neighborhood.

	\subsection{\label{sec:movement}The vehicles movement.} 	In order to~define vehicle traffic rules and simulate their movement, the~model proposed by~\citeauthor{LarAlvIca2010:Cellular}~(\citeyear{LarAlvIca2010:Cellular}) was used. The~proposed model meets the~general behavior of~vehicles on~the~road. Drivers with free space ahead are traveling at~maximum speed. Approaching the~second vehicle, drivers react to~changes in~its speed, providing themselves with a~constant space for collision-free braking. This model is~often called \textbf{LAI~model}, from the~authors' names. This part of~the~work will include a~description of~this model and also comments on~possible assumptions.
	
	The model presents traffic flow at~a~single-lane road, where vehicles move from left to~right. The~road is~divided into 2.5-meters sections, and each is~presented as~a~separate cell. The~length of~the~car is~taken as~5 meters what is~represented as~two cells. Each cell can be~empty or~occupied only by~part of~one vehicle. The~position of~the~vehicle is~determined by~the~position of~its front bumper. Vehicles run at~speeds from 0 to~$v_{max}$, which symbolize the~number of~cells a~vehicle can move in~one-time step $t$. The~time step corresponds to~one second. The~speed conversion from simulation to~real is~presented in~the~Table \ref{tab:speeds_simulations}. 
	\begin{table}[h!]
	\caption{The relationship between real and simulation speeds in~the~model used.}
	\label{tab:speeds_simulations}
	\centering
		\begin{tabular}{|c|c|c|c|} \hline
			\textbf{Velocity $v$} & \textbf{Distance} & \textbf{Real speed}  & \textbf{Real speed}\\ \hline
			1 & 2.5 m & 2.5 m/s & 9 km/h \\  \hline
			2 & 5 m & 5 m/s & 18 km/h \\ \hline
			3 & 7.5 m & 7.5 m/s & 27 km/h \\  \hline
			4 & 10 m & 10 m/s & 36 km/h \\  \hline
			5 & 12.5 m & 12.5 m/s & 45 km/h \\  \hline
			6 & 15 m & 15 m/s & 54 km/h \\  \hline
			7 & 17.5 m & 17.5 m/s & 63 km/h \\  \hline
		\end{tabular}
	\end{table}
	
	Here in~the~first column, we have the~velocity used in~the~model, next column presents how distance is~done in~a~one-time step (1 second), the~next columns present real velocity in~m/s and km/h for better imagine how the~model works. In~simulations we decide to~used maximum speed equals to~45 km/h, because traffic flow in~the~city is~considered, so~drivers have not too much space to~fast driving. 
	
	The model takes into account the~limited acceleration and braking capabilities of~vehicles and also ensures appropriate distances between vehicles to~guarantee safe driving. Three distances calculated for the~car following to~its predecessor are included. These values calculate the~distance needed for safe driving in~the~event that the~driver wants to~slow down ($d_{dec}$), accelerate ($d_{acc}$) or~maintain the~current speed ($d_{keep}$), assuming that the~predecessor will want to~suddenly start slowing down with the~maximum force M until to~stop. They~are calculated as~follows:

\begin{subequations}\label{Ddist}{\scriptsize
	\begin{align}
d_{acc} &= \max \Big( 0,\sum^{\floor{(v_n(t)+\Delta v)/M}}_{i = 0} \hspace{-2em}\left[ (v_n(t) + \Delta v) - i\cdot M\right] - \sum^{\floor{(v_{n+1}(t)-M)/M}}_{i = 0}  \hspace{-2em}\left[ (v_{n+1}(t) - M) - i\cdot M \right]  \Big) \label{dist1} \\ &\nonumber \vspace{-4ex} \\
d_{keep} &= \max \Big( 0, \sum^{\floor{v_n(t)/M}}_{i = 0} \left[ v_n(t) - i\cdot M\right] - \sum^{\floor{(v_{n+1}(t)-M)/M}}_{i = 0}  \left[ (v_{n+1}(t) - M) - i\cdot M \right]  \Big) \label{dist2} \\ 
&\nonumber \vspace{-4ex} \\
d_{dec} &= \max \Big( 0, \sum^{\floor{(v_n(t)-\Delta v)/M}}_{i = 0} \hspace{-2em}\left[ (v_n(t) - \Delta v) - i\cdot M\right]- \sum^{\floor{(v_{n+1}(t)-M)/M}}_{i = 0}  \hspace{-2em}\left[ (v_{n+1}(t) - M) - i\cdot M \right]  \Big) \label{dist3}
	\end{align}}
\end{subequations}

Here, vehicle $n$ is~the~follower, and $n + 1$ is~the~preceding car. $v_n(t)$ means the~value of~the~velocity of~vehicle $n$ in~time $t$, $\Delta v$ is~the~ability to~accelerate in~one-time step and $M$ is~ability to~emergency braking. 
	
	Updating vehicle traffic takes place in~four steps, which are done parallel for each of~the~vehicles. 
	\begin{enumerate}[I.]\itemsep1pt \parskip1pt \parsep1pt
		\item Calculation of~safe distances $d_{dec_n}$, $d_{acc_n}$, $d_{keep_n}$. 
		\item  Calculation of~the~probability of~slow acceleration.
		\item Speed update.
		\item Updating position.
	\end{enumerate}
	
	\textbf{Safe distances.} According to~formulas \ref{Ddist} safe distances are counted for each vehicles. The~calculation of~these values is~based on~the~assumption that if~the~vehicle in~the~next time step $t + 1$ increases its speed (or maintains it~or~ slows it~down respectively) and the~driver preceding from the~moment $t$ will constantly slow down to~speed 0 (with maximum ability to~emergency braking), there will be~no collision. The~different between these equation is~just in~first part, which define traveled distance by~vehicle $n$ if~it~decelerate ($v_n(t+1) = v_n(t) - \Delta v$), keep velocity ($v_n(t+1) = v_n(t)$) or~accelerate $v_n(t+1) = v_n(t) + \Delta v$, in~next time step, and next begins to~brake rapidly. The~second part of~equation determines the~distance traveled by~the~preceding vehicle if~it~starts to~braking with maximum force~$M$. 
	\vspace{0.8ex}
	
	\textbf{Calculation of~the~probability of~slow acceleration.} 
	The second step in~the~vehicle movement procedure focuses on~calculating the~stochastic parameter $R_a$ responsible for slowing down vehicle acceleration. It is~assumed that low-speed vehicles have more troubles to~accelerate. According to~human nature and the~mechanism of~the~car, it~is~true that is~that the~faster we go, the~easier we manage to~accelerate, and standing or~driving very slowly cause slower acceleration. The~limiting speed at~which acceleration comes easier is~assumed to~be 3, which corresponds to~27 km/h. The~value of~$R_a$ parameter is~calculated based on~the~formula
	\begin{equation}
	R_a = \min (R_d, R_0 + v_n(t) \cdot (R_d - R_0)/v_s),
	\end{equation}
	where $R_0$ and $R_d$ are fixed stochastic parameters, mean respectively probability to~accelerate when the~speed is~equal to~0, and probability to~accelerate when the~speed is~equal or~more than $v_s$, and $v_s$ limit speed below which acceleration is~harder. 
	
\begin{figure}[tbh!]
		\centering
		\includegraphics[height=4cm,width=9cm]{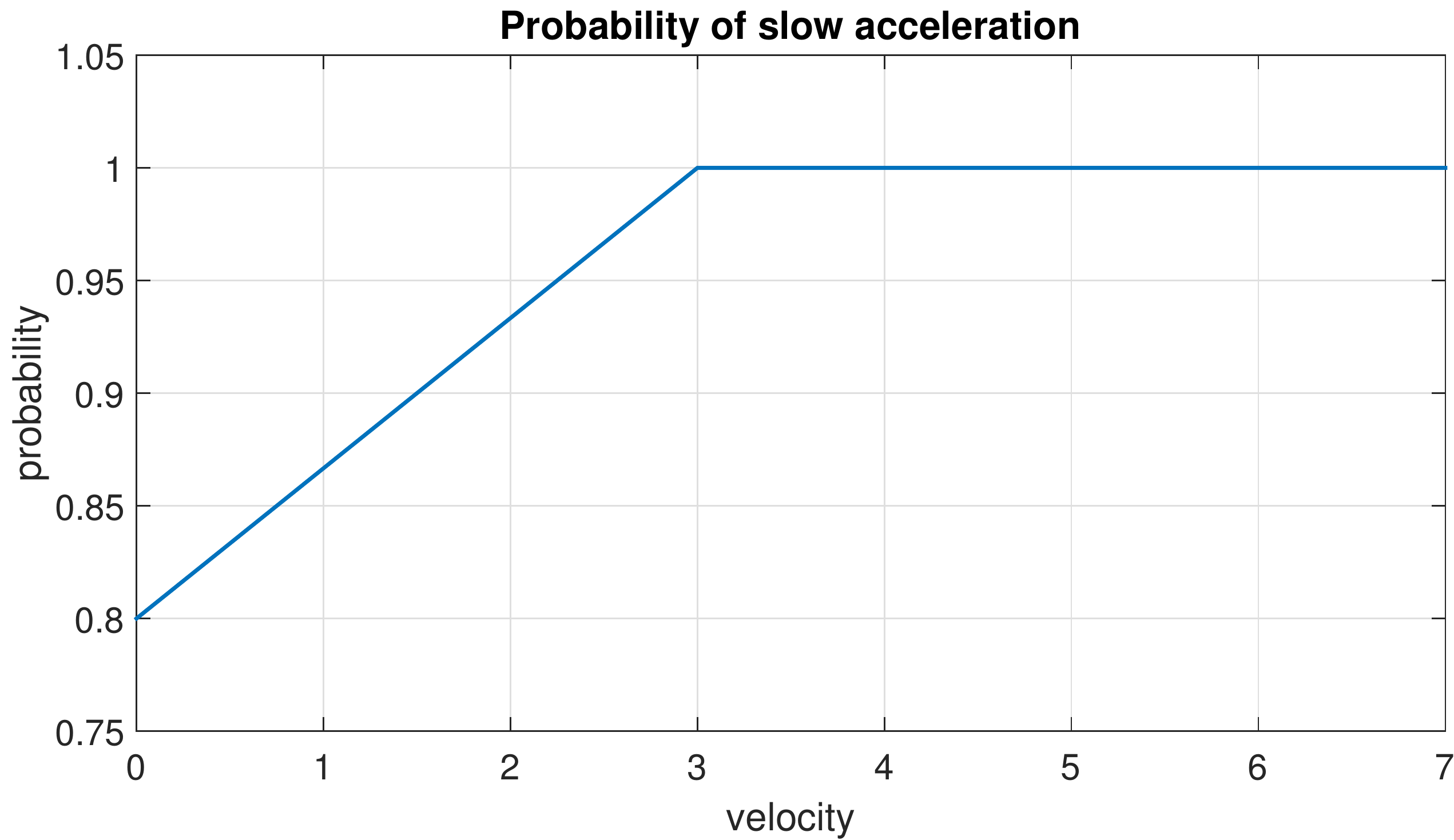} 
		\caption{Values of~$R_a$ parameters for fixed $R_0$, $R_d$ and $v_s$ \cite{LarAlvIca2010:Cellular}.} \label{Ra_lai_plot}
	\end{figure}

	Easily can be~seen,  that the~relationship between the~probability of~acceleration at~speed 0 and at~a~speed greater than the~limit is~interpolated linearly, which is~taken from the~idea presented also by~\citeauthor{Ra_ref} \cite{Ra_ref}. In~the~simulations, 0.8 and 1 were adopted as~$R_0$ and $R_d$ parameters, respectively, which will not cause frequent difficulties in~accelerating vehicles, however, the~stochastic nature of~this process will be~taken into account. The~graph of~the~$R_a$ parameter change for the~other parameters adopted in~this way is~presented in~Figure \ref{Ra_lai_plot}.  
	
	\vspace{0.8ex}
	\textbf{Speed update.}
	In the~beginning, as~mentioned before $\Delta v$ means speed increase in~one time step, fixed for all vehicles. $v_n(t)$ and $x_n(t)$ determine the~velocity and the~position of~vehicle $n$ in~time $t$. Distance from vehicle $n$ to~vehicle $n+1$ is~counted by~the~following formula 
	\begin{equation*}
	d_n(t) = x_{n+1}(t) - x_n(t) - l_s,
	\end{equation*}
	which exactly means the~distance from front bumper of~vehicle $n$ pointed by~$x_n(t)$ to~rear bumper of~the~vehicle in~the~front, presented by~the~difference between the~position of~the~front bumper $x_{n+1}(t)$ and the~length of~the~vehicle $l_s$ (in~cells).
	The speed update is~done in~four steps, the~order of~which does not matter.
	\begin{enumerate}\itemsep1pt \parskip1pt \parsep1pt
		\item \textit{Acceleration.} If the~distance to~the~preceding vehicle is~greater than $d_{aac_n}$ then the~vehicle $n$ increase velocity by~$\Delta v$ with probability $R_a$, what is~presented as~follows
		\begin{equation*}
		v_n (t + 1) = \left\{
		\begin{array} {ll}
		\min (v_n (t) + \Delta v,v_{max}), & \textrm{with prob. $ R_a $} \\
		v_n (t), & \textrm {otherwise} \\
		\end{array} \right.
		\end{equation*}
		In this rule is~assumed that all drivers strive to~achieve the~maximum velocity if~it~is~possible. Here is~include irregular ability to~accelerate depends on~the~distance to~preceding vehicles, relevant velocities of~both, and stochastic parameter responsible for slower acceleration defined in~Step II.  
		\item \textit{Random slowing down.} This rule allows drivers to~maintain the~current speed, if~it~allows safe driving, it~also takes into account traffic disturbances, which are an indispensable element of~traffic flow. The~probability of~random events is~determined by~the~$R_s$ parameter. If $d_{acc_n} > d_n(t) \geq d_{keep_n} $, then the~updated speed is~determined according to~the~formula 
		\begin{equation*}
		v_n (t + 1) = \left\{
		\begin{array} {ll}
		\max (v_n (t) - \Delta v,0), & \textrm{with prob. $ R_s $} \\
		v_n (t), & \textrm {otherwise} \\
		\end{array} \right.
		\end{equation*}
		\item \textit{Braking.} This rule ensures that the~drivers keep an adequate distance from front vehicles. Rapid braking is~not desirable, so~in~order to~ensure a~moderate braking process for the~driver, when the~free space in~front of~the~car is~too small, the~vehicle speed is~reduced by~$\Delta v$, which reflects optimal braking.    
		\begin{equation*}
		v_n(t+1) = \max(v_n(t) - \Delta v,0) \quad \textrm{if } \quad d_{keep_n} > d_n(t) \geq d_{dec_n}
		\end{equation*}
		\item \textit{Emergency braking.} As can be~seen in~real life, it~is~not always possible to~brake calmly. What is~more, road situations often force more aggressive braking.  Such situations are included in~this rule. When the~driver gets too close to~the~other car, or~when the~other car brakes too much, it~forces emergency braking. If the~distance is~at least $d_{dec}$, this rule is~not applied. According to~the~commonly accepted standard proposed in~the~literature (v.  \citeauthor{force_M} \citeauthor{LarAlvIca2010:Cellular}, the~emergency braking force is~set to~$-5 $ m/s$^2$. With respect to~the~assumed model parameters the~value of~$M$ is~2. This step is~described by~ 
		\begin{equation*}
		v_n(t+1) = \max(v_n(t) - M,0) \quad \textrm{if } \quad  d_n(t) < d_{dec_n}
		\end{equation*}
	\end{enumerate}
	
	\vspace{0.8ex}
	\textbf{Updating position} Finally, with updated vehicle speed, it~is~possible to~actualize vehicle positions. The~vehicles are moved by~the~number of~cells according to~their speed. This is~described by~means of~
	\begin{equation*}
	x_n(t+1) = x_n(t) + v_n(t+1),
	\end{equation*}
	where $x_n(t+1)$ is~actualized position, $v_n(t+1)$ is~the~previously determined vehicle speed, and $x_n(t)$ is~last position of~vehicle.

	\subsection{Intersections}\label{sec:intersections}
	Intersections are an inseparable element of~road traffic, they are an intersection with a~road at~one level. All connections and crossroads also count in~intersections. There are the~following types of~intersections: 
	\begin{itemize}\itemsep1pt \parskip1pt \parsep1pt
		\item uncontrolled intersections
		\item intersections with traffic signs
		\item crossings with controlled traffic (traffic lights or~authorized person)
	\end{itemize}
	Modeling of~traffic at~intersections is~an important element of~road traffic modeling, many models have been created on~this subject, such as~models simulating the~movement of~vehicles at~intersections of~type T \cite{T_shape1}, describing the~movement at~un-signalized intersections as~in~the~case of~\cite{intersection1}, \cite{intersections2} and those considering traffic at~intersections with traffic lights \cite{signalized1}. Typically, these models consist of~two aspects, modeling vehicle traffic and modeling interactions at~intersections. General rules are set for intersections, however, the~behavior of~drivers who may or~may not comply with these rules is~also taken into account. Modeling of~such behavior is~also different, which usually distinguishes these models. This aspects was consider in~my engineering thesis 
	\cite{SzaWlo2020:Divers}.  Helpful in~modeling interactions at~intersections is~game theory, which facilitates the~decision about the~right of~way, where players are drivers in~conflict at~the~intersection, examples of~such use can be~seen in~\cite{gamet1}. Additionally, signalized intersection models using Markov chain are often used, as~in~the~case of~\cite{markov1}.

	However, the~purpose of~the~work is~simple modeling of~road traffic, therefore advanced intersection modeling methods will not be~considered. It is~assumed that all drivers comply with traffic regulations and follow road safety. The~consequences of~changing behavior to~incorrect and inconsistent with traffic rules are not investigated. The~purpose of~the~work is~to~find elements that affect the~potential threat affecting the~reluctance of~drivers to~comply with traffic rules. Depending on~the~maneuver performed by~the~drivers and the~type of~intersection, the~following situations need to~be modeled:
	\begin{itemize}\itemsep1pt \parskip1pt \parsep1pt
		\item turn right from the~road without right of~way
		\item turn right from the~road with right of~way
		\item turn left from the~road with right of~way
		\item turn left from the~road without right of~way 
		\item go straight ahead at~traffic lights
		\item turn left at~traffic lights
	\end{itemize}
	When modeling the~above situations, two basic rules were used:
	\begin{description}\itemsep1pt \parskip1pt \parsep1pt
		\item[Rule 1] a~driver who wants to~join the~traffic on~the~main road can do, if~and only if, during the~whole process, until the~maximum speed is~reached, he does not disturb the~driving of~other vehicles on~the~main road. This maneuver may be~described by~the~following formula: 
		\begin{equation*}
		l_x - v_x - \sum_{\Delta v=2}^{v_{max}}\min(v_{max},v_x+\Delta v-1) + \sum_{\Delta v=2}^{v_{max}} \Delta v > d_{keep_x},
		\end{equation*}
		where $l_x$ is~the~distance of~the~vehicle on~the~main road to~the~intersection, $v_x$ is~his current vehicle speed. The~first sum symbolizes the~distance traveled by~the~vehicle on~the~main road until the~passing vehicle reaches maximum speed. The~second sum represents the~distance traveled by~the~vehicle joining the~traffic until it~reaches maximum speed, assuming that in~the~first second the~vehicle will be~at an intersection with a~speed equal 1. Both vehicles increase their speed by~1 in~each second and they do not exceed the~maximum speed. The~value of~the~left side of~the~inequality must be~greater than the~distance needed by~the~driver on~the~main road to~maintain his speed. Otherwise, the~driver would be~forced to~brake which would disturb his movement.
		\item[Rule 2] The~driver wanting to~cross the~opposite direction road can do it~if there is~no collision with the~opposite direction during the~time needed to~complete it~and the~opposite driver will not be~forced to~brake. The~time it~takes to~complete the~maneuver depends on~the~initial speed at~the~start of~the~maneuver. This~relationship is~described in~the~Table \ref{tab:time_cross}.
		\begin{table}[H]
		    \centering
		    \caption{Relationship between the~time of~crossing of~the~opposite road and the~initial speed.}
		    \label{tab:time_cross}
		    \begin{tabular}{|c|c|}
				\hline 
				Velocity $v_n$ & Need time $\tau_n$ \\ 
				\hline 
				2 & 1s \\ 
				\hline 
				1 & 2s \\ 
				\hline 
				0 & 3s \\ 
				\hline 
			\end{tabular} 
		\end{table}
		The condition ensuring the~correct execution of~the~maneuver can be~described by~the~following inequality:
		\begin{equation*}
		l_x - \sum_{\Delta v=0}^{\tau_n}\min(v_{max},v_x+\Delta v-1) > d_{dec_x},
		\end{equation*}
		where $l_x$ is~the~distance of~the~opposite car opposite to~the~intersection and the~sum is~responsible for calculating the~distance traveled by~this car in~the~time needed to~complete the~turn. The~value of~the~left side of~the~inequality must be~greater than the~speed needed for safe braking by~the~vehicle. Otherwise, it~would force the~driver to~emergency braking, which is~not desirable, and in~the~event of~a~possible delayed reaction of~the~driver could lead to~an accident.
	\end{description}
	The above rules are the~basis used to~define behavior at~intersections, more information on~where the~rules were applied, and why, it~will be~described in~Section \ref{sec:street_describtion}.
	
\begin{figure}[tbh!]
		\centering
		\includegraphics[height=3cm]{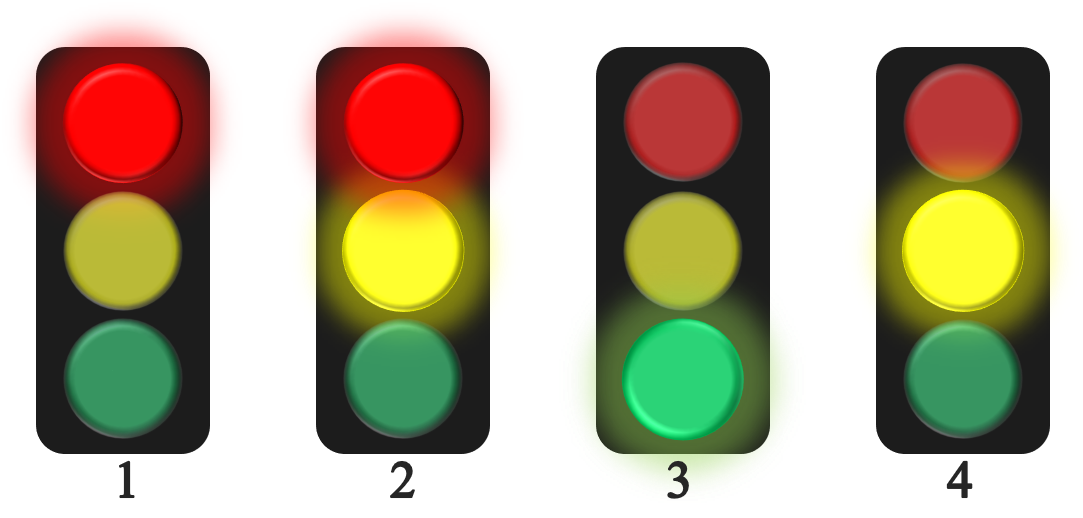} 
		\caption{Traffic light cycle.} \label{sygnalizator}
	\end{figure}

	In addition, modelling of~traffic lights was needed. The~traffic light scheme was used in~accordance with the~Polish regulations. The~traffic light cycle follows the~diagram \ref{sygnalizator}.
	The duration and meaning of~individual signals are as~follows:
	\begin{enumerate}\itemsep1pt \parskip1pt \parsep1pt
		\item \textbf{Red light} --- no entry behind the~signal light. The~duration is~60 seconds.
		\item \textbf{Red and yellow light} --- means that in~a~moment will be~a~green signal.  According to~the~regulations, it~lasts 1 s. 
		\item \textbf{Green light} --- allows entry after the~signal light if~it~is~possible to~continue driving and this will not cause a~road safety hazard. The~duration is~the~same as~for the~red signal, equal to~60 s.
		\item \textbf{Yellow light} --- does not allow entry behind the~signal light, unless stopping the~vehicle would cause an emergency brake. According to~the~regulations, it~should last at~least 3 seconds. 
	\end{enumerate}
	Such a~traffic light cycle and the~duration of~each signal were adopted in~the~simulation. of~course, there is~also a~relationship between the~capacity of~intersections and the~time of~the~traffic light cycle. However, the~most standard signaling scheme was adopted to~ensure optimal intersection capacity. In~addition, it~was assumed that both directions of~travel are equivalent, which is~why this cycle is~the~same on~both roads.

	In accordance with the~theory described for traffic modeling, as~well as~with the~proposed method of~traffic conditioning at~intersections. For each street from the~diagram in~the~drawing \ref{fig:main_structure}, traffic simulations were performed in~the~MATLAB package. Real and simulated street sizes are presented in~Section \ref{sec:street_describtion} in~the~next chapter. For each simulation, the~time it~took me from the~beginning of~the~road to~leaving the~intersection at~its end was calculated for each vehicle. Simulations have been carried out many times for different probabilities of~a~new driver appearing on~the~road, which in~the~further understanding will be~taken as~traffic intensity. 
	
	\subsection{Model calibration}
	An important aspect in~the~case of~traffic modeling, which we could not fail to~mention in~this chapter, is~the~calibration of~models. In~general, this topic is~part of~a~larger problem, which is~simulation optimization. The~area of~development of~simulation optimization in~recent years has enjoyed great interest among researchers and practitioners. Simulation optimization is~the~pursuit of~the~maximum performance of~a~simulated real system. The~system performance is~assessed based on~the~simulation results, and the~model parameters are the~decision variables. The~assessed performance in~this case is~the~model's ability to~recreate reality. Therefore, it~is~a~very important topic in~modeling traffic, which aims to~enable the~reconstruction of~real vehicle traffic, so~correctly choose the~model parameters so~that the~model used is~a~reliable model and correctly shows the~modeled behavior characteristics. Optimization in~the~context of~motion simulation models has evolved in~many areas and the~only ones were optimization and calibration of~motion, but they were not often combined with optimization theory, where some of~the~problems in~motion modeling are well known. One of~the~most important conclusions is~that there is~no algorithm that is~suitable for all problems and needs and that the~choice of~the~right algorithm depends on~the~example being examined (v. \citeauthor{Spall2006}(\citeyear{Spall2006})).
Most studies focused on~testing the~performance of~the~optimization algorithm, where models are evaluated against actual traffic data, e.g. \Citeauthor{Hollander2008} in~\citeyear{Hollander2008}. However, based on~real traffic data, it~is~not possible to~evaluate the~effectiveness of~the~algorithm and the~entire calibration process. Another approach proposed in~the~literature is~the~use of~synthetic measurements, i.e. data obtained from the~model itself. This approach was proposed e.g. by~ \citeauthor{Ossen2008}(\citeyear{Ossen2008}), and tested changes in~model calibration due to~the~use of~errors in~synthetic motion trajectories by~\citeauthor{ciuffo}(\citeyear{ciuffo}), which used tests with synthetic data to~configure the~process of~calibration of~microscopic motion models, based on~trial and error.

	\section{\label{chap:simulations}Simulation} 
	\subsection{\label{sec:street_describtion}Description of~real traffic network} 
	Using the~models proposed in~section~\ref{chap:trafic_modelling}, a~simulation of~vehicle movement was performed on~each street presented in~the~Figure \ref{fig:main_structure}. The~model of~vehicle traffic along a~straight road is~presented in~Section \ref{sec:movement}. The~modeling movement between streets was more complicated. Section \ref{sec:intersections} describes the~general rules needed to~define traffic at~intersections. There are, various maneuvers required simulation. In~addition, the~actual road lengths have been converted into simulation values to~best reflect the~road traffic. Table \ref{opis_ulic_tabela} describes real and simulation road lengths and maneuvers that should be~performed on~a~given road section.
	\begin{table}[tbh!]
		\centering{\small
		\caption{Description and base information on~analysed roads.}
		\label{opis_ulic_tabela}
		\begin{tabular}{|c|c|c|c|p{6cm}|} \hline
			\multirow{2}{0.5cm}{\textbf{Id}} & \textbf{Street} &
			\multicolumn{2}{c|}{\textbf{Length}} &
			\multirow{2}{7cm}{\textbf{~~~~~Intersections and turning}} \\
			\cline{3-4}
			& \textbf{name} & \textbf{In meters} & \textbf{In cells} & \\ \hline
			1 & Dolna & 300 m & 120 & Give way on~Zlota and turn right\\  \hline
			2 & Zlota & 350 m & 140 & Give way oncoming vehicles and turn left or~go straight \\ \hline
			3 & Mickiewicza & 500 m & 200 & Turn right with right of~way\\  \hline
			4 & Piwna~1 & 450 m & 180 & Give way oncoming vehicles and turn left or~go straight\\  \hline
			5 & Nyska~1 & 160 m & 64 & Go ahead with right of~way \\  \hline
			6 & Zlotnickiego & 500 m & 200 & Give way on~Nyska~1 and turn right \\ \hline
			7 & Piwna~2 & 180 m & 72 &  Give way vehicles on~the~main road and turn left\\  \hline
			8 & Jasna & 400 m & 160 &  and Give way vehicles on~the~main road and turn left \\  \hline
			9 & Nyska~2 & 200 m & 80 & Give way oncoming vehicles and turn left on~intersection \\  \hline
			10 & Laska & 500 m & 200 & Go ahead, but wait on~traffic lights\\  \hline
			11 & Sieradzka~1 & 500 m & 200 & Go ahead with right of~way\\  \hline
			12 & Sieradzka~2 & 500 m & 200 & Go ahead to~the~end of~road\\  \hline
		\end{tabular}}
	\end{table}
	In the~case of~Nyska~1, Sieradzka~1, and Sieradzka~2 streets, drivers go through given section with priority, driving straight ahead.
	
	For Dolna and Zlotnickiego roads, drivers join the~traffic on~the~main road being on~a~road without the~right of~way. Rule 1 was applied, assuming that when approaching an intersection, drivers must slow down to~a~speed of~0 or~1, and then decide according to~the~condition described.
	
	At Mickiewicza street, at~the~end of~the~road, the~driver is~forced to~slow down to~2, which corresponds to~the~real speed of~18 km/h, we can assume that this is~a~reasonable speed to~make a~turn. After decelerating, drivers can leave the~intersection.
	
	For Zlota and Piwna~1 streets, drivers with probability $\frac{1}{3} $ turn left, otherwise they go straight. Before turning, the~drivers slow down to~at least speed 2, if~they can cross the~opposite direction lane they continue driving if~they do not slow down more. Therefore, drivers must give way to~oncoming vehicles, Rule 2 applies.
	
	In the~case of~Piwna~2 and Jasna streets, we assume that the~drivers slow down before the~intersection to~0 or~1 and with probability $\frac{1}{2} $ turn right or~left. In~both situations, it~is~necessary to~apply Rule 1, because drivers must give way to~vehicles that are on~the~road they are turning, in~addition in~the~case of~a~left turn, Rule 2 should be~applied too because the~vehicle will cross the~opposite direction.
	
	The last two traffic situations relate to~traffic at~intersections with traffic lights. When driving along Laska Street, in~the~event of~red light, drivers wait at~the~intersection, then they can leave it. The~case where drivers would like to~turn left is~not being considered because in~real life a~left lane is~intended for a~left turn.  When leaving Nyska~2 Street, drivers may ride to~the~right, left, or~straight. In~the~case of~a~right turn or~straight ahead, the~process goes without any problems, so~we allow drivers to~leave the~intersection. When turning left, you must pass vehicles driving in~the~opposite direction, so~Rule 2 applies.
	
	According to~the~above assumptions, simulations were made, and repeated 1000 times for each traffic intensity to~obtain the~average values of~delays depending on~the~intensity. The sample code and description of the program are at the end of the work in Appendix~\ref{chap:appendixCode}. 

	\subsection{Simulations results}
\begin{figure}[h!]
		\centering
		\includegraphics[height=5cm,width=11cm]{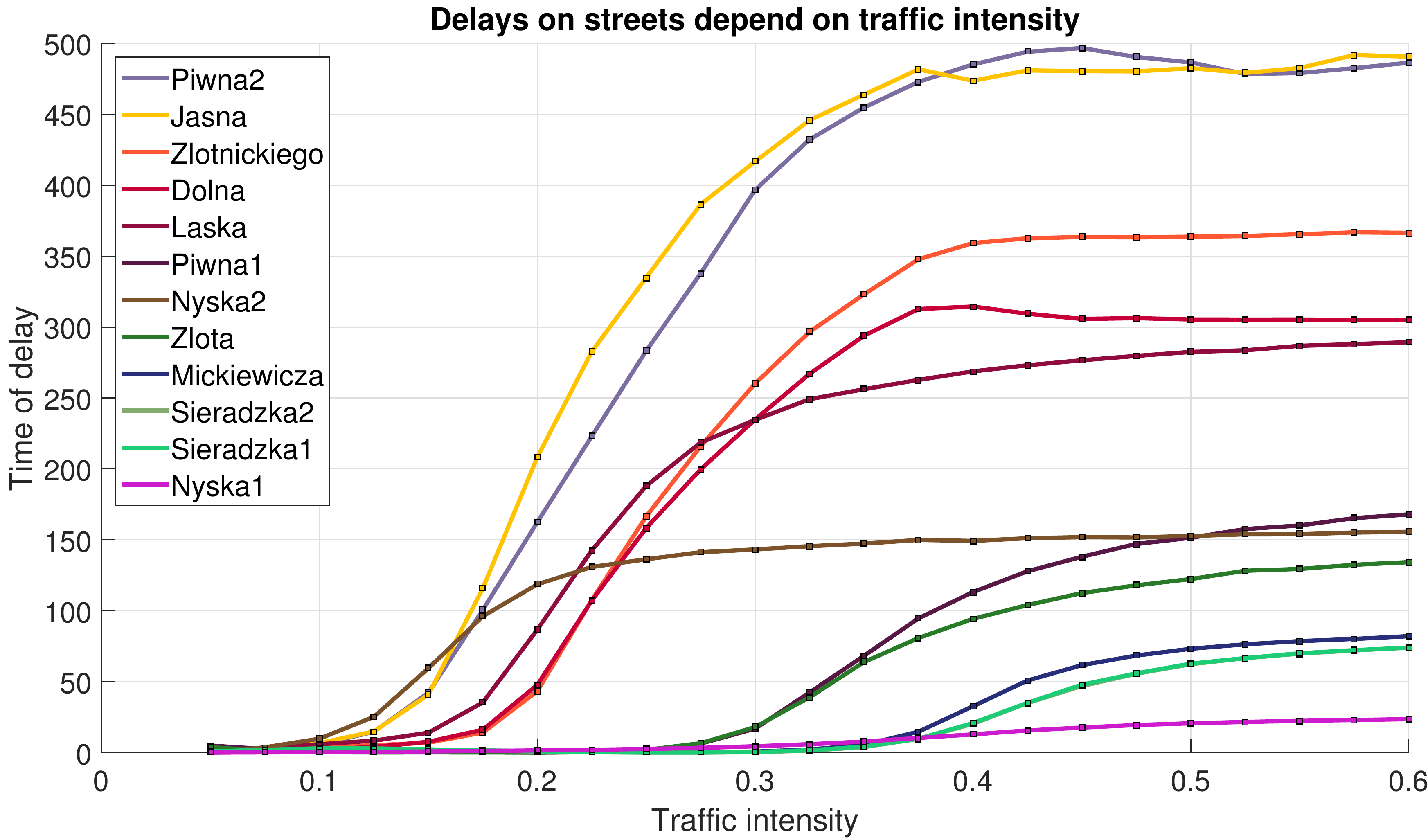} 
		\caption{\label{delays}Time of~delays depend on~traffic intensity for different roads.} 
	\end{figure}

	In accordance with the~characteristics described in~the~previous section, simulations of~motion were made. The~results of~road delays are shown in~graph \ref{delays}. 
	We see that the~delay increase characteristics for different roads are different. It is~easy to~see that one of~the~most difficult streets to~travel are Jasna and Piwna~2, we see here a~high sensitivity to~traffic intensity. Another group of~streets in~terms of~delays are Zlotnickiego and Dolna, and also Laska Street is~similar to~them, although the~growth characteristics are different. Nyska~2 has a~completely different behavior from the~rest, but it~is~the~only street with such a~complex intersection, including traffic lights. In~this case, the~delay increases very quickly, reaching a~critical level for this street, related to~the~capacity of~the~road. Therefore, despite the~fact that the~final result of~the~delay is~not the~largest, it~can be~considered that the~efficiency of~this intersection is~the~worst. The~next, but definitely more efficient streets are Piwna~1 and Zlota, and the~most fluid traffic can be~seen on~the~last 4 streets, where there are no intersections and traffic disturbances. In addition, both Sieradzka streets have the same delay times, because they both are without intersections streets and they have the same length.
	
\begin{figure}[tbh!]
		\centering
		\includegraphics[height=5cm,width=11cm]{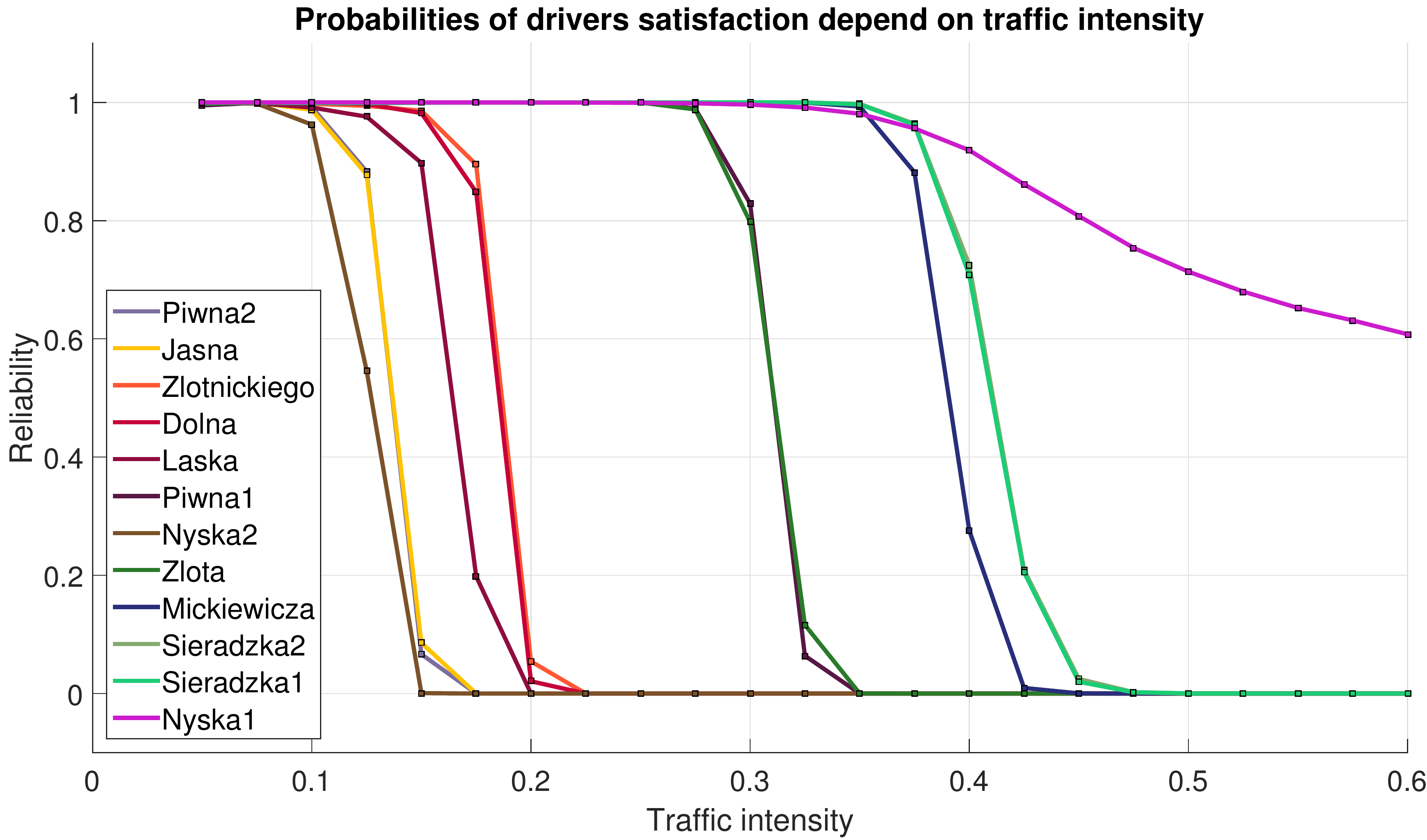} 
		\caption{\label{reliabilities_plot}Probabilities of~drivers' satisfaction vs.~traffic intensity for different roads.} 
	\end{figure}

	Next, using the~approach introduced in~Section \ref{sec:reliabilityImp_traffic} and based on~the~calculated delay times, it~is~possible to~determine the~probability of~driver satisfaction with a~given section of~the~route. An undesirable phenomenon is~exceeding a~certain critical level of~delay time, which will cause dissatisfaction to~the~driver. The~probability that the~critical value for a~given delay is~not exceeded is~described by~the~reliability function of~the~Weibull distribution. Using this, the~probability of~driver satisfaction for a~given delay on~each road will be~calculated depending on~traffic intensity. These probabilities are presented in~Figure \ref{reliabilities_plot}.
	\hspace{-2ex}
	We see that the~reliability of~individual streets is~different. Ones of~the~streets already at~a~low traffic intensity reach a~critical state, which will cause drivers' dissatisfaction for sure~(these are e.g. ulice Nyska~2, Piwna~2, Jasna, Laska, Dolna, Zlotnickiego). We~can also observe streets such as~Nyska~1, where traffic is~constantly flowing and does not irritate drivers. This drawing shows us that it~is~true that individual streets react differently to~increasing traffic. Therefore, it~is~worth examining how these changes affect the~overall functioning of~the~traffic network and the~importance of~individual fragments.
	
	\subsection{Values of~Importance Measures}
	To begin with, we calculate the~reliability of~the~entire system depending on~traffic and the~reliability of~each route. In~this way, we obtain the~probability of~finish the~journey with the~satisfaction of~all roads. In~order to~calculate the~reliability values of~individual roads $ p_i $, for $ i = 1,2, \ldots, 12$ are substituted the~appropriate values in~the~formula \eqref{eq:my_struc_fun}. In~addition, we will calculate the~reliability of~individual routes that correspond to~the~minimum paths. The~relationship between the~elements of~each route is~in~series. Individual routes include the~following streets:
	\begin{description}\itemsep1pt \parskip1pt \parsep1pt
		\item[Route 1:] Dolna, Zlota, Mickiewicza, Jasna, Sieradzka~2
		\item[Route 2:] Dolna, Zlota, Nyska~1, Nyska~2, Sieradzka~1, Sieradzka~2
		\item[Route 3:] Piwna~1, Zlotnickiego, Nyska~2, Sieradzka~1, Sieradzka~2
		\item[Route 4:] Piwna~1, Piwna~2, Laska, Sieradzka~1, Sieradzka~2
	\end{description}
	
	\vspace{-2ex}
	\begin{table}[h!]
		\centering{\small
		\caption{	\label{niezawodnosci}The probabilities of~comfortable driving between points $A$ and $B$ vs. various routes and different traffic intensity.}
	\hspace{-1ex}
		\begin{tabular}{|c|c|c|c|c|c|} \hline
			\textbf{Traffic} & \multicolumn{5}{c|}{\textbf{Probability of~satisfaction from}} \\
			\cline{2-6}
			\textbf{intensity} & \textbf{All routes} & \textbf{Route 1} & \textbf{Route 2} & \textbf{Route 3} & \textbf{Route 4}\\ \hline
			\textbf{0.050} & 1.0000 & 0.9974  & 0.9973 & 0.9948  & 0.9948 \\ \hline
			\textbf{0.075} & 1.0000 &  0.9989 & 0.9970  & 0.9961  &  0.9969\\ \hline
			\textbf{0.100} & 1.0000 & 0.9855  & 0.9598  & 0.9574  & 0.9807 \\ \hline
			\textbf{0.125} & 0.9963 & 0.8737  & 0.5440  & 0.5422   & 0.8596 \\ \hline
			\textbf{0.150} & 0.1396 & 0.0841  &  0.0006 & 0.0006  & 0.0595 \\ \hline
			\textbf{0.175} & 0 & 0  & 0  &  0  &  0 \\ \hline
			\textbf{0.200} & 0 & 0 & 0 & 0 & 0 \\ 
			$\vdots$ & $\vdots$ & $\vdots$ & $\vdots$ & $\vdots$ & $\vdots$ \\ 
			\textbf{0.600} & 0 & 0 & 0 & 0 & 0 \\ \hline
		\end{tabular}}
	\end{table}
	The calculated reliability values are presented in~Table \ref{niezawodnosci}. We can see that the~system is~no longer efficient at~a~traffic intensity of~0.175. In~addition, we see that the~capacity of~the~system is~always greater than the~efficiency of~individual roads. This is~important information regarding the~critical value of~traffic intensity that causes failure of~the~entire network. In~addition, we can see that the~capacity of~Routes 1 and 4 is~greater, which may suggest that with heavy traffic it~is~better to~choose one of~these two routes to~ensure a~better chance of~a~quiet ride. \\

	We will now proceed to~calculate the~importance of~individual roads in~the~functioning of~the~entire system. The~calculated values are shown in~Table \ref{importances_values}. For each street, received values of~measure of~significance at~a~given traffic intensity were presented. As mentioned before, these values are calculated on~the~basis of~the~structure function \eqref{eq:my_struc_fun} and importance measures theory introduced by~Birnbaum \eqref{eq:Birnbaum_Imp1} using the~received reliability for individual traffic intensities.
	\begin{table}[tbh!]
		\centering{\small
		\caption{\label{importances_values}The importance of~route elements for different traffic intensities.}
		\hspace{1ex}
		\begin{tabular}{|c|c|c|c|c|c|c|c|c|c|} \hline
			\multirow{2}{0.5cm}{\textbf{Id}} & \textbf{Street} & \multicolumn{8}{c|}{\textbf{Traffic intensity}} \\ 
			\cline{3-10}
			& \textbf{name} & \textbf{0.050} &  \textbf{0.075}  &  \textbf{0.100}  &  \textbf{0.125}  &  \textbf{0.150}  & \textbf{0.175}  &  \textbf{0.200} & $\ldots$ \\ \hline
			1 & Dolna & $\approx 0$ & $\approx 0$ & $\approx 0$ & 0.0301 & 0.0810 &  0& 0 & $\ldots$ \\  \hline
			2 & Zlota & $\approx 0$ & $\approx 0$ & $\approx 0$ & 0.0300 & 0.0795 &  0& 0 & $\ldots$ \\  \hline
			3 & Mickiewicza & $\approx 0$ & $\approx 0$ & $\approx 0$ & 0.0256 & 0.0790 & 0 & 0 & $\ldots$ \\  \hline
			4 & Piwna~1 & $\approx 0$ & $\approx 0$ & $\approx 0$ & 0.0271 & 0.0550 & 0 & 0 & $\ldots$ \\  \hline
			5 & Nyska~1 & $\approx 0$ & $\approx 0$ & $\approx 0$ & 0.0044 & 0.0005 & 0 & 0 & $\ldots$ \\  \hline
			6 & Zlotnickiego & $\approx 0$ & $\approx 0$ & $\approx 0$ & 0.0044 & 0.0005 & 0 & 0 & $\ldots$ \\  \hline
			7 & Piwna~2  & $\approx 0$ & $\approx 0$ & $\approx 0$ & 0.0257 & 0.8201 & 0.1980 & 0 & $\ldots$ \\  \hline
			8 & Jasna & $\approx 0$ & $\approx 0$ & $\approx 0$ & 0.0292 & 0.9214 & 0.8481 & 0.0207 & $0,\ldots$ \\  \hline
			9 & Nyska~2 & $\approx 0$ & $\approx 0$ & $\approx 0$ & 0.0161 & 1.6921 & 1.7424 & 0.0751 & $0,\ldots$ \\  \hline
			10 & Laska & $\approx 0$ & $\approx 0$ & $\approx 0$ & 0.0232 & 0.0607 & 0 & 0 & $\ldots$ \\  \hline
			11 & Sieradzka~1 & $\approx 0$ & $\approx 0$ & $\approx 0$ & 0.0315 & 0.0556 & 0 & 0 & $\ldots$ \\  \hline
			12 & Sieradzka~2 & $\approx 0$ & $\approx 0$ & 0.0001 & 0.0571 & 0.1346 & 0 & 0 & $\ldots$ \\  \hline
		\end{tabular}}
	\end{table}
	As we can see, the~most interesting results were obtained for the~traffic intensity of~0.125 and 0.150. At low traffic intensities, the~reliability of~individual elements does not affect the~functioning of~the~system, because the~whole system works properly and the~reliability of~the~roads are close to~1. For the~intensity of~0.125, the~contribution of~individual streets begins to~be noticeable. We see that, according to~structural measures, the~largest contribution to~the~functioning of~the~network has Sieradzka~2 street, the~next streets have the~value of~importance close to~0.03, with the~exception of~Nyska~1, Zlotnickiego, and Nyska~2, which are smaller. For the~intensity of~0.150, we can see that there are difficulties in~movement. The~first thing that draws our attention is~the~importance of~Nyska~2 Street, which was one of~the~smallest before, now it~has become the~most significant element. Another important component of~the~system is~again Sieradzka~2, which is~obvious. However, the~streets that are worth paying attention to~are Piwna~2 and Jasna, whose significance has also risen dramatically. With subsequent increases in~intensity, we see that only these 3 streets really affect the~quality of~traffic, and of~them the~most street Nyska~2.

	\subsection{Comparison with real life data.}
	Based on~the~analysis made in~the~previous section, the~streets Nyska~2, Piwna~2, and Jasna have the~greatest importance for the~appropriate functioning of~the~entire system at~high traffic intensity. The~analyzed traffic system is~a~real traffic network, which is~why we know what traffic really looks like on~individual roads. The~presented scheme of~travel from A~to B shows the~travel from two strategic positions in~the~city. The~main streets in~the~city are Laska and Sieradzka, they pass through the~center of~the~city. The~traffic "on top" of~Laska Street is~greater than on~Sieradzka Street because here we are already approaching the~exit from the~city. The~results obtained are in~line with expectations. One of~the~most important points in~the~city is~the~Nyska--Laska--Sieradzka intersection. In~fact, this intersection is~more extensive and we can see that a~lot of~work has been put into its proper functioning. Many simplifications are used there, which would also slightly change the~results obtained from the~simulation. For example, vehicles turning left into Sieradzka Street have more space so, when they are waiting for a~turn, they do not obstruct the~traffic of~other vehicles going straight or~turning right. In~addition, time counters are used on~the~traffic lights that increase drivers' watchfulness and their start when the~green light comes on.
	
	
	The intersection of~Piwna and Laska streets was critical enough that it~was impossible to~turn left there, the~sign 'right to~turn right' was in~force. This was a~major impediment to~general traffic as~well as~to~the~routes presented in~the~paper. That is~why a~roundabout intersection has recently been built here. This decision certainly required a~lot of~consideration by~the~city authorities, because there is~not enough space for a~full-size roundabout here, so~it~has a~slightly flattened one side. However, as~can be~seen from the~results obtained, it~was one of~the~critical parts of~traffic in~the~city, so~this decision seems sensible.
	
	The last problematic street is~Jasna, but here in~real traffic, there is~no such intensity of~vehicles, both on~Jasna Street and the~"bottom" part of~Sieradzka Street. Assuming that traffic in~this part of~the~city is~smaller and that turning into Mickiewicza street is~not very problematic gives important information to~drivers who considered which of~these two roads is~better.

	\section{Summarizing and conclusions.}	A~quantitative approach to~road quality assessment is proposed. The measures of~significance defined for the reliability systems were used as~a~tool to~calculate the~importance of~individual road fragments. An~actual traffic network diagram was analyzed, ensuring access from point $A$ to~point $B$. To begin with, assuming that only the~structure of~the~analyzed road network is~known, the~structural significance of~individual road fragments was calculated. For this purpose, two measures were used: proposed by~Birnbaum, assuming constant reliability of~individual road fragments, and Barlow and Proschan measure, which takes into account the~variability of~individual element reliability. There were noticeable differences between the~received values, but the~final result was similar in~both cases. The~most important for maintaining the~efficiency of~the~analyzed road network are the~streets that occur in~the~largest number of~possible routes to~the~point $ B $ as~a~serial connection. This~confirmed our expectations, but also helped to~locate some of~the~roads that at~the~first consideration were not potentially important routes. When comparing the~differences between Birnbaum and Barlow-Proschan measures, the~second one was considered more appropriate for use in~the~context of~road traffic, because the~reliability of~individual roads are not the~same, many factors affect on~them. 
	
	Then a~method of~assessing the~reliability of~street elements was proposed. For this purpose, it~was assumed that the~quality of~roads is~the~satisfaction of~drivers with the~route traveled, and the~delay time on~individual roads was used as~a~measure of~this. It was assumed that drivers have limited patience, which is~close to~a~lifetime and was presented as~a~variable from the~Weibull distribution. Having calculated the~delay times on~individual roads, it~was possible to~determine the~probability of~upset the~driver at~such a~delay. However, in~order for the~obtained value to~be able to~be used in~the~theory of~measures of~importance, it~had to~be transformed so~that it~was responsible for the~reliability of~a~given element. Therefore, the~Weibull distribution reliability function was used, which in~our example reflected the~probability that with a~given road delay, the~driver would still be~satisfied. Alternate method of~reliability assessing to~the~net of~roads has been used by~\citeauthor{PilSzy2009:Road}(\citeyear{PilSzy2009:Road}).


In the paper by \citeauthor{SzaWlo2020:Divers}~(\citeyear{SzaWlo2020:Divers})), it~was shown that if~drivers are dissatisfied with driving then they can stop complying with traffic rules. And one of~the~factors influencing their change and negative behavior on~the~road is~delays. Therefore, the~measures defined in~this way are a~guide for both drivers and traffic managers. For drivers, it~shows which road is~better to~avoid because there is~a~chance of~potential nervousness, and gives road drivers information about dangerous points in~the~city and points that have a~negative impact on~drivers. In~addition, the~large delay time on~individual roads indicates the~failure of~the~fragments concerned. Based on~the~simulations performed, the~delay times on~each road were calculated. Then it~was shown which road fragments are the~most important. The~obtained results were confronted with the~actual feelings regarding the~given fragments. And they were considered likely because with the~network as~defined it~was used as~the~most significant elements that were improved in~real traffic. Which proves the~real importance of~these elements. 


\vspace{6pt} 
\authorcontributions{\noindent \textbf{Author Contributions:} Development of the application of the important measures to element of road networks, Krzysztof Szajowski(KSz) and Kinga Włodarczyk(KW); implementation of the algorithm and the example, numerical simulations, KW; writing  and editing, KSz and KW.
}

\funding{Supported by Wrocław University of Science and Technology, Faculty of Pure and Applied Mathematics, under the project 049U/0051/19(KSz).
}


\conflictsofinterest{The authors declare no conflict of interest.
} 

\abbreviations{The following abbreviations are used in this article:\\

\noindent 
\begin{tabular}{@{}ll}
\textbf{TDL} & Time-Dependent Lifetime (p. \pageref{TDL})\\
\textbf{TIL} & Time Independent Lifetime (p. \pageref{TDL})\\
\textbf{LAI model}& Model of  general behavior of vehicles on the road (p. \pageref{sec:movement}; v. \cite{LarAlvIca2010:Cellular})
\end{tabular}}


	
	\nopagebreak

	\bibliography{\jobname}
	\bibliographystyle{abbrvnat}
	
\appendix

\lstset{language=Matlab,%
		basicstyle=\scriptsize,
		breaklines=true,%
		morekeywords={matlab2tikz},
		keywordstyle=\color{blue},%
		morekeywords=[2]{1}, keywordstyle=[2]{\color{black}},
		identifierstyle=\color{black},%
		stringstyle=\color{mylilas},
		commentstyle=\color{mygreen},%
		showstringspaces=false,
		numbers=left,%
		numberstyle={\tiny \color{black}},
		numbersep=9pt, 
		emph=[1]{for,end,break},emphstyle=[1]\color{red}, 
	}
	\section{Code of modelling and simulation\label{chap:appendixCode}}
	This work uses the LAI model presented in Section \ref{sec:movement}. Which was then modified to add restrictions at intersections. Each of the streets included in the intersection was simulated individually and the interactions at the intersection were examined. As it was mentioned before, streets and intersections can be divided into several types. In this appendix, we will present the most advanced in terms of regulations used, i.e. the intersection on Piwna ~2 and Jasna streets. Cars there may go to the right joining the traffic in this direction or turn left crossing the second direction of travel in addition. Both rules apply at intersections. Code presentations will start with the most basic ones and then we will go to the main code. To calculate the distance to the vehicles preceding, which will provide the possibility of acceleration, maintain speed or deceleration was calculated using three functions \verb|d_acc.m|, \verb|d_keep.m|, \verb|d_dec.m|, listed below
\begin{lstlisting}
function result = d_acc(n,speeds,M,dv)
% n - id of car
% speeds - vector of speeds
% M - max ability to decelerate
% dv - ability do accelerate

result = max(0,sum((speeds(n) + dv) - (0:(floor((speeds(n)+dv)/M)))*M) - ...
    sum((speeds(n+1) - M) - (0:(floor((speeds(n+1)-M)/M)))*M));

end

\end{lstlisting}
\begin{lstlisting}
function result = d_keep(n,speeds,M,dv)

result = max(0,sum(speeds(n) - (0:(floor(speeds(n)/M)))*M) - ...
    sum((speeds(n+1) - M) - (0:(floor((speeds(n+1)-M)/M)))*M));

end
\end{lstlisting}
\begin{lstlisting}
function result = d_dec(n,speeds,M,dv)

result = max(0,sum((speeds(n) - dv) - (0:(floor((speeds(n)-dv)/M)))*M) - ...
    sum((speeds(n+1) - M) - (0:(floor((speeds(n+1)-M)/M)))*M));

end
\end{lstlisting}
	Another important element is adding a new vehicle on the road, each vehicle has its own index, speed, position, and information where it goes. Depending on the simulation being performed, the probability of route selection can be set to others, if it is not needed, the \textit{where} parameter is not given.
\begin{lstlisting}
function [pos, speeds, ids, where] = 
               new_car(pos,speeds,ids,v_min,prob_new_car,M,dv,ls,L,where)
% where - where the vehicle are going - 1 right, 0 - left
if nargin < 10 
    where = [];
end
% Generating a new car, with some probability. New cars are added with
% v_min velocity
if isempty(pos)
    if rand() <= prob_new_car
        pos = [2,pos];
        speeds = [v_min,speeds];
        where = [rand()>1/2, where];
        ids = 1;
    end
else
    temp = 1:(L/2);
    if pos(1) ~= 2
        speeds_temp = [v_min,speeds];
        d_ke = d_keep(1,speeds_temp,M,dv); 
        % space sufficient to maintain current speed
        if pos(1) - ls >= d_ke 
        % if the distance to the previous one is sufficient 
        % to maintain the current speed then you can enter
            if rand() <= prob_new_car
                pos = [2,pos];
                speeds = [v_min,speeds];
                where = [rand()>1/2, where]; % 1 right, 0 left
                ids = [min(temp(~ismember(1:(L/2),ids))),ids];
            end
        end
    end
end
end


\end{lstlisting}
	Then, a speed update is performed for each vehicle according to the diagram described in chapter \ref{sec:movement}. In addition, parameters are used to say whether the vehicle can leave the intersection or not, and to what speed it should slow down before the intersection.
\begin{lstlisting}
function [speeds,can_go] = velocity_updade(n, speeds, pos, ls, dv, Rd, R0, Rs
                                          , vs, M, vmax, L, can_go, dec_to)
% Velocity update for vehicle 'n'
% can_go --- 
% 0 - the car cannot leave the intersection, 
% 1 - the car can leave the intersection
% dec_to --- velocity do decelerate before intersection

if n ~= numel(pos) % for last car no dependence of previous car
    d_ac = d_acc(n,speeds,M,dv);
    d_ke = d_keep(n,speeds,M,dv);
    d_de = d_dec(n,speeds,M,dv);
    % Acceleration
    if (pos(n+1) - pos(n) - ls) >= d_ac
        new_speeds = speeds(n);
        if rand()<= min(Rd,R0+speeds(n)*(Rd-R0)/vs)
            new_speeds = min(speeds(n)+dv,vmax);
        end
        speeds(n) = new_speeds;
    end
 % Random slowing down  
 if d_ac > (pos(n+1) - pos(n) - ls) &&  (pos(n+1) - pos(n) - ls) >= d_ke
        new_speeds = speeds(n);
        if rand() <= Rs
            new_speeds = max(speeds(n)-dv,0);
        end
        speeds(n) = new_speeds;
    end
% Braking
 if d_ke > (pos(n+1) - pos(n) - ls) && (pos(n+1) - pos(n) - ls) >= d_de
        new_speeds = max(speeds(n)-dv,0);
        speeds(n) = new_speeds;
    end
    % Emergency braking
    if (pos(n+1) - pos(n) - ls) < d_de
        new_speeds = max(speeds(n)-M,0);
        speeds(n) = new_speeds;
    end
else
    speeds_temp = [speeds,dec_to];
    d_ke = d_keep(n,speeds_temp,M,dv);
    d_de = d_dec(n,speeds_temp,M,dv);
    if d_ke < (L - pos(n)) || can_go
        % Acceleration
        new_speeds = speeds(n);
        if rand()<= min(Rd,R0+speeds(n)*(Rd-R0)/vs)
            new_speeds = min(speeds(n)+dv,vmax);
        end
        speeds(n) = new_speeds;
        % Random slowing down
        new_speeds = speeds(n);
        if rand() <= Rs
            new_speeds = max(speeds(n)-dv,0);
        end
        speeds(n) = new_speeds;
        can_go = 0; % swich of for next car
    end
    % Slow before end od road
    if d_ke > (L - pos(n)) && (L - pos(n)) >= d_de
        new_speeds = max(speeds(n)-dv,0);
        speeds(n) = new_speeds;
    end
    % Emergency braking
    if (L - pos(n)) < d_de
        new_speeds = max(speeds(n)-M,0);
        speeds(n) = new_speeds;
    end
end
end
\end{lstlisting}
	Finally, we go to the main program codes. Depending on the exact type of intersection, the program looks slightly different, but the overall characteristics and construction are preserved. At the beginning we define the variables used in the model, then we have loops after repetitions for different probabilities of a new vehicle. We define new empty roads in each loop. Then we add the first vehicle on each road and go on to further processes. In the original program, before starting the loop after repeating the update on the road, the road was filled with vehicles. At each step, we update speeds and add a new vehicle on the road, according to the probability. Then update the speeds and remove those vehicles whose position has exceeded the length of the road. Finally, we analyze the interactions of drivers at intersections, resulting in a change in the parameter saying whether the vehicle can leave the intersection or not. This is done in accordance with the previously described assumptions. On the posted program we have an example for Piwna Street~2, where the driver turning right gives way to other vehicles on this road and turning left gives way to vehicles driving in the opposite direction than he plans because he crosses their lane. Below is the code.
\begin{lstlisting}
MC = 1000; % Monte Carlo
dx = 2.5; % cell size
ls = 2; %  car length, 5 m
vmax = 5; % 5 -> 12.5 m/s -> 45 km/h lub 6 -> 15 m/s -> 54 km/h
Rs = 0.01; % prob. of emergency
Rd = 1; % max. prob of acceleration
R0 = 0.8; % min. prob. of acceleration
vs = 5/dx + 1; % min velocity for faster driving
M = 2; % ability to emergency driving 5m/s2
dv = 2.5/dx; % ability to accelerate
dec_to = 0; % velocity to slow down before intersection

same_probs = 0.05:0.025:0.6; % different traffic intesity
it = 1;

for same_prob = same_probs
 disp(same_prob)
 SAVE = zeros(1,MC);
 for M = 1:MC
 % Intersection Piwna2 i Laska, z prawd. 1/2 turn right or left
 % Inputs for Piwna2
 L_Piwna2 = 180/dx; % road length
 prob_new_car_Piwna2 = same_prob;
 pos_Piwna2 = [];
 speeds_Piwna2 = [];
 v_min_Piwna2 = 4;
 can_go_Piwna2 = 0;
       
 % Inputs for Laska1, from right
 L_Laska1 = 100/dx;
 prob_new_car_Laska1 = same_prob;
 pos_Laska1 = [];
 speeds_Laska1 = [];
 v_min_Laska1 = 4;
 can_go_Laska1 = 1;
        
 % Inputs for Laska2, from left
 L_Laska2 = 100/dx;
 prob_new_car_Laska2 = same_prob;
 pos_Laska2 = [];
 speeds_Laska2 = [];
 v_min_Laska2 = 4;
 can_go_Laska2 = 1;
        
 % First car on Piwna2, Laska1 and Laska2
 if isempty(pos_Piwna2) % if on road there is no car
            pos_Piwna2(1) = 2;
            speeds_Piwna2(1) = v_min_Piwna2;
            ids_Piwna2 = 1;
            where_Piwna2 = [1]; % 1 - right, 0 - left
 end
        
 if isempty(pos_Laska1) % if on road there is no car
            pos_Laska1(1) = 2;
            speeds_Laska1(1) = v_min_Laska1;
            ids_Laska1 = 1;
 end
        
 if isempty(pos_Laska2) % if on road there is no car
            pos_Laska2(1) = 2;
            speeds_Laska2(1) = v_min_Laska2;
            ids_Laska2 = 1;
 end
        
 point = L_Laska1; 
 % intersection point on Laska1 i Laska 2, at the end od roads
        stop = zeros(1,3);
        licz_czas_Piwna2 = zeros(size(1:(L_Piwna2/2)));
        licz = 1;
        saved_czas_Piwna2 = [];

        can_go_Piwna2 = 0; 
 % must always stop before the intersection
 for i = 1:1000
 % Adding new car with v_min velocity on Piwna2
 [pos_Piwna2, speeds_Piwna2, ids_Piwna2, where_Piwna2] = new_car(pos_Piwna2
    ,speeds_Piwna2,ids_Piwna2,v_min_Piwna2,prob_new_car_Piwna2,M,dv,ls
    ,L_Piwna2, where_Piwna2);
 % Velocity updating for each vehicle on road Piwna2
 for n = 1:numel(pos_Piwna2)
    [speeds_Piwna2, ~ ] = velocity_updade(n, speeds_Piwna2, pos_Piwna2, ls
    , dv, Rd, R0, Rs, vs, M, vmax, L_Piwna2,can_go_Piwna2,dec_to);
 end
            
 % Adding new car with v_min velocity on Laska1
    [pos_Laska1, speeds_Laska1, ids_Laska1] = new_car(pos_Laska1
    ,speeds_Laska1,ids_Laska1,v_min_Laska1,prob_new_car_Laska1,M
    ,dv,ls,L_Laska1);
% Velocity updating for each vehicle on road Laska1
    for n = 1:numel(pos_Laska1)
    speeds_Laska1 = velocity_updade(n, speeds_Laska1, pos_Laska1, ls
    , dv, Rd, R0, Rs, vs, M, vmax, L_Laska1,can_go_Laska1,dec_to);
    end
            
% Adding new car with v_min velocity on Laska2
    [pos_Laska2, speeds_Laska2, ids_Laska2] = new_car(pos_Laska2
    ,speeds_Laska2,ids_Laska2,v_min_Laska2,prob_new_car_Laska2,M
    ,dv,ls,L_Laska2);
% Velocity updating for each vehicle on road Laska2
    for n = 1:numel(pos_Laska2)
    speeds_Laska2 = velocity_updade(n, speeds_Laska2, pos_Laska2
    , ls, dv, Rd, R0, Rs, vs, M, vmax, L_Laska2,can_go_Laska2,dec_to);
    end
            
% Position update
    pos_Piwna2 = pos_Piwna2 + speeds_Piwna2;
    pos_Laska1 = pos_Laska1 + speeds_Laska1;
    pos_Laska2 = pos_Laska2 + speeds_Laska2;
            
            
% Removing cars
    if sum(pos_Piwna2 > L_Piwna2)
       speeds_Piwna2(pos_Piwna2 > L_Piwna2) = [];
       where_Piwna2(pos_Piwna2 > L_Piwna2) = [];
       ids_Piwna2(pos_Piwna2 > L_Piwna2) = [];
       saved_czas_Piwna2(licz) = max(licz_czas_Piwna2);
       licz = licz + 1;
       licz_czas_Piwna2(licz_czas_Piwna2 == ...
                     max(licz_czas_Piwna2)) = 0;
       pos_Piwna2(pos_Piwna2 > L_Piwna2) = [];
     end
            
     if sum(pos_Laska1 > L_Laska1)
                speeds_Laska1(pos_Laska1 > L_Laska1) = [];
                pos_Laska1(pos_Laska1 > L_Laska1) = [];
                
     end
     if sum(pos_Laska2 > L_Laska2)
                speeds_Laska2(pos_Laska2 > L_Laska2) = [];
                pos_Laska2(pos_Laska2 > L_Laska2) = [];
     end
licz_czas_Piwna2(ids_Piwna2) = licz_czas_Piwna2(ids_Piwna2) + 1;
% Intersection interaction
 if ~isempty(speeds_Piwna2)
    if sum(pos_Piwna2 == L_Piwna2) 
% some vehicle is at the end of the road
% Different conditions for cars which go in right and in left
    if where_Piwna2(end) == 1 % in right
     to_point_on_left = pos_Laska2 - point;
     finds_on_left = find(to_point_on_left<0);
     if ~isempty(finds_on_left)
     found_on_left = find(to_point_on_left(finds_on_left) == ...
            max(to_point_on_left(finds_on_left))); 
% founded car on Zlota is the closest to intersections
    dist_to_point_on_left = abs(pos_Laska2(found_on_left) - point);
    d_keep_on_left = max(0,sum(vmax - (0:(floor(vmax/M)))*M) - ...
    sum((vmax - M) - (0:(floor((vmax-M)/M)))*M)); 
% distance which will not force the driver 
% to be released on the main road at the maximum speeds of both vehicles
    if dist_to_point_on_left - speeds_Laska2(found_on_left) - ...
    sum(min(vmax,speeds_Laska2(found_on_left)+(2:vmax)-1)) + ... 
    sum(2:vmax) >= d_keep_on_left
% if the condition is met it can enter
   can_go_Piwna2 = 1; 
% the vehicle before the intersection can already go
        else
    can_go_Piwna2 = 0;
  end
  else
  can_go_Piwna2 = 1;
  end

else % where=0, in left
cond = 0;
% Checks if it is free on the right
to_point_on_right = pos_Laska1 - point;
finds_on_right = find(to_point_on_right<0);
if ~isempty(finds_on_right)
    found_on_right = find(to_point_on_right(finds_on_right) == ...
    max(to_point_on_right(finds_on_right))); 
    % founded car on Zlota is the closest to intersections
    dist_to_point_on_right = abs(pos_Laska1(found_on_right) - point);
    d_keep_on_right = max(0,sum(vmax - (0:(floor(vmax/M)))*M) - ...
    sum((vmax - M) - (0:(floor((vmax-M)/M)))*M)); 
% distance which will not force the driver 
% to be released on the main road at the maximum speeds of both vehicles
 if sum(pos_Piwna2 == L_Piwna2) % some vehicle is at the end of the road
                                
 if dist_to_point_on_right - ... 
       sum(min(vmax,speeds_Laska1(found_on_right)+(0:2))) - ... 
       sum(min(vmax,speeds_Laska1(found_on_right)+(2:vmax)-1)) + ...
       sum(2:vmax) >= d_keep_on_right
% if this condition is met on the right we have 
%  the opportunity to enter the intersection
       cond = cond + 1;
    end
            end
    else
    cond = cond + 1;
           end
                        
% Checking if it is free on the left
to_point_on_left = pos_Laska2 - point;
finds_on_left = find(to_point_on_left<0);
if ~isempty(finds_on_left)
found_on_left = find(to_point_on_left(finds_on_left) == ... 
max(to_point_on_left(finds_on_left)));
dist_to_point_on_left = abs(pos_Laska2(found_on_left) - point);
    speeds_temp = [speeds_Laska2(found_on_left),0];
    d_ke_on_left = d_keep(1,speeds_temp,M,dv);
    d_de_on_left = d_dec(1,speeds_temp,M,dv);
    needed_time = 3;
if dist_to_point_on_left - speeds_Laska2(found_on_left) - ... 
    min(speeds_Laska2(found_on_left)+1,vmax) - ...
    min(speeds_Laska2(found_on_left)+2,vmax) >= d_de_on_left
    cond = cond + 1; % second condition is met
end
else
cond = cond + 1; % second condition is met
end
                        
if cond == 2 % if both conditions were met vehicle can go
   can_go_Piwna2 = 1;
   else
   can_go_Piwna2 = 0;
                        end
                    end
                end
            end
        end
SAVE(M) = nanmean(saved_czas_Piwna2(5:end));
    end
    SAVE_MC(it) = nanmean(SAVE(SAVE>0));
    it = it+1;
end
\end{lstlisting}
	Similarly, all simulations carried out at work were carried out, thus obtaining the travel time for each section. On this basis, the road delay was calculated as the difference between each reading, the smallest value.

\end{document}